\documentclass[moor]{informs1}              

\usepackage{natbib,bm,bbm}
\usepackage[framemethod=tikz]{mdframed}
\usepackage{graphicx}\usepackage{subfig,tikz,pgfplots,float,makecell}\usetikzlibrary{arrows}
\usepackage{amsmath}

 \NatBibNumeric
 \def\bibfont{\small}%
 \bibpunct[, ]{[}{]}{,}{n}{}{,}%

\usepackage[colorlinks=true,breaklinks=true,bookmarks=true,urlcolor=blue,
     citecolor=blue,linkcolor=blue,bookmarksopen=false,draft=false]{hyperref}

\def\EMAIL#1{\href{mailto:#1}{#1}}
\def\URL#1{\href{#1}{#1}}         

\TheoremsNumberedThrough     

\EquationsNumberedThrough    

\def\proof{\noindent{\bf Proof}\ \ }
\def\qed{\hfill \vrule height 7pt width 7pt depth 0pt\medskip}
\newcommand{\ds}{\displaystyle}

\newcommand{\ba}{\begin{array}}
\newcommand{\ea}{\end{array}}

\newcommand{\1}{\mathbbm{1}}
\renewcommand{\l}{\left}\renewcommand{\r}{\right}
\renewcommand{\be}{\begin{equation}}
\renewcommand{\ee}{\end{equation}}
\newcommand{\eps}{\varepsilon}

\renewcommand{\mc}{\mathcal}

\newcommand{\ov}{\overline}
\newcommand{\ul}{\underline}

\newcommand{\R}{\mathbb{R}}

\newcommand{\C}{\mathbb{C}}

\newcommand{\se}{\text{ if }}

\def\R{\mathbb{R}}
\def\C{\mathbb{C}}

\begin{document}


 \RUNAUTHOR{Massai, Como, and Fagnani}

\RUNTITLE{Equilibria and Systemic Risk in Saturated Networks}

 \TITLE{Equilibria and Systemic Risk in Saturated Networks}

\ARTICLEAUTHORS{%
\AUTHOR{Leonardo Massai, Giacomo Como,\thanks{Giacomo Como is also with the Department of Automatic Control, Lund University, BOX 118, SE-22100, Lund, Sweden.} and Fabio Fagnani  }
\AFF{Dipartimento di Scienze Matematiche, Politecnico di Torino, Corso Duca degli Abruzzi 24, 10129, Torino, Italy \EMAIL{leonardo.massai@polito.it}, \EMAIL{giacomo.como@polito.it}, \EMAIL{fabio.fagnani@polito.it}\\ 
\URL{https://staff.polito.it/giacomo.como/},  \URL{https://staff.polito.it/fabio.fagnani/}}
} 

\ABSTRACT{%
We undertake a fundamental study of network equilibria modeled as solutions of fixed point  equations for monotone linear functions with saturation nonlinearities.
The considered model extends one originally proposed to study systemic risk in networks of financial institutions interconnected by mutual obligations 
 and is one of the simplest continuous models accounting for shock propagation phenomena and cascading failure effects. 
It also characterizes Nash equilibria of constrained quadratic network games with strategic complementarities. 
 We first derive explicit expressions for network equilibria and prove necessary and sufficient conditions for their uniqueness encompassing and generalizing  results available in the literature. Then, we study jump discontinuities of the network equilibria when the exogenous flows cross certain regions of measure $0$ representable as graphs of continuous functions. 
Finally, we discuss some implications of our results in the two main motivating applications. In financial networks, this bifurcation phenomenon is responsible for how small shocks in the assets of  a few nodes can trigger major aggregate losses to the system and cause the default of several agents. In constrained quadratic network games, it  induces a blow-up behavior of the sensitivity of Nash equilibria with respect to the individual benefits.
}%


\MSCCLASS{}
\ORMSCLASS{Primary: ; secondary: }
\HISTORY{}

\maketitle

%


%
%
%

\section{Introduction} \label{sec:intro}


A central aspect of complex socio-technical systems such as infrastructural, social, economic, and financial networks is the role played by interconnections in amplifying and propagating shocks through cascading mechanisms that may increase the fragility of a system  \cite{Allen.Gale:2000,Buldyrev.ea:2010,Como.ea:2013b}. The term \emph{systemic risk} refers to the possibility that even small shocks localized in a limited part of the network can spill over, thus possibly achieving a significant global impact \cite{Haldane.May:2011,Savla.ea:2014,Acemoglu2015a}. 
A key challenge is to find adequate models for network systems, that are sufficiently elaborate to incorporate such propagation phenomena, yet simple enough to allow for mathematical tractability. Whilst simple contagion models such as epidemic contact processes prove inadequate as they are based on purely pairwise interactions, more complex models taking into account cumulative neighborhood effects include the linear threshold model \cite{Adler:91,Watts:2002,Rossi.ea:2019} whose applicability is however limited by the fact that states of the nodes are described by pure binary variables simply expressing whether the node has been affected by the shock. 

In most of the applications  where the network represents a physical infrastructure, a social or economic network, or an interconnected financial system, however, the cascading mechanism is rather triggered by a process naturally described in terms of continuous variables such as, e.g., power flows in electric grids, traffic volumes in transportation systems, the extent of individuals' involvement  in a certain activity in social networks, prices or quantities of goods in an economic system, assets values and payments in financial networks. 
The most tractable continuous models of network interaction considered in the literature give rise to notions of equilibria that can be mathematically characterized as the solutions of a linear system of equations whose coefficients can be assembled in a (often sparse) matrix that describes the network of interconnections among the different nodes. Examples include competitive equilibria in production networks \cite{Acemoglu.ea:2012,Acemoglu2015} or Nash equilibria  in network games with linear best replies including quadratic network games \cite{Ballester.ea:2006,CalvoArmengol.ea:2009,Galeotti.ea:2010,Bramoulle.ea:2014,Jackson.Zenou:2015}. 

While the most basic formulations of  such fundamental models consider no constraints on the involved variables, in several of the aforementioned applications it is natural to assume some a priori lower (e.g., non-negativity) and upper bounds (e.g., maximum available resource). 
E.g., in the Eisenberg and Noe model \cite{Eisenberg.Noe:2001}, financial institutions are interconnected by mutual obligations and the payments are necessarily non-negative and upper bounded by the debt value. In the context of network games modeling peer effects on students' engagement,  
\cite{CalvoArmengol.ea:2009} suggests to ``bound the strategy space in such a game rather naturally by simply acknowledging the fact that students have a time constraint and allocate their time between leisure and school work,'' and 
\cite{Bramoulle.Kranton:2016} acknowledge that 
``while in principle, a player's action could be any real number, all games in the literature place restrictions on players' actions which represent different real-world situations'' and that ``for peer effects in a classroom, there are natural lower and upper bounds: a student can study no less than zero hours and no more than twenty-four hours in a day.''
When a priori upper and lower bounds are taken into account in the network model, the related equilibrium notions end up being mathematically characterized as the solutions of linear systems of equations with saturation non-linearities.  
\cite{CalvoArmengol.ea:2009,Belhaj.ea:2014,Bramoulle.ea:2014,Allouch:2015,Bramoulle.Kranton:2016}. Such saturated network models exhibit a considerably richer behaviors than purely linear ones, including the possibility of cascading effects coded in terms of variable saturations and transition phases with respect to structure parameters.

%
In this paper, we undertake a fundamental study of such saturated equilibrium models in networks with positive externalities. 
Precisely, we consider the following fixed point equation
\begin{equation}\label{eq1}
x_i=\min \left\{ \max \left\{ \sum_{j=1}^n x_jP_{ji}+c_i, 0 \right\} , w_i \right\}\,,\qquad i=1,\ldots,n\,,
\end{equation}
where $P$ in $\R_+^{n\times n}$ is a non-negative square matrix  and $w$ in $\R_+^n$ is a non-negative vector  that jointly describe a \emph{network}, while $c$ in $\R^n$ is an \emph{exogenous flow} vector.
Equation \eqref{eq1} can be more compactly rewritten as
\begin{equation}\label{eq:nonlinear-system}
x=S_0^w\left(P^{\top}x+c\right)\,.
\end{equation}
where $S_0^w$ denotes the vector saturation function
\begin{equation}\label{eqd}
\left(S_0^w(x)\right)_i = \min \{ \max \{ x_i, 0\}, w_i\}, \qquad i=1,\ldots,n\,,
\end{equation}
We shall refer to vectors $x$ that are solutions of  \eqref{eq:nonlinear-system} as \emph{equilibria} of the network $(P,w)$ with exogenous flow $c$. 
Notice that the range of the vector saturation function $S_0^w$ is contained in the complete lattice
\be\label{lattice}\mc L_0^w=\{x\in\R^n:\,0\le x\le w\}\,.\ee
As the lattice $\mc L_0^w$ is a nonempty, convex, and compact set, and $x\mapsto S_0^w(P^{\top}x+c)$ maps $\mc L_0^w$ in itself with continuity,
existence of network equilibria directly follows from Brower's fixed point theorem. Hence, the set $\mc X\subseteq \mc L_0^w$ of network equilibria is always nonempty. On the other hand, the structure of such network equilibria as well as their uniqueness and dependence on the exogenous flow prove to be more delicate issues. They will be the object of this paper.

In financial networks, starting with the seminal work of Eisenberg and Noe \cite{Eisenberg.Noe:2001},  the entries of the vector $w$ represent the obligations of the various institutions, those of the exogenous flow $c$ represent the balance between assets possessed by the entities and their obligations towards institutions outside the network, while $P$ is a row-stochastic or sub-stochastic matrix describing the way obligations of an entity are split among the others thus encoding the backbone of the financial system interconnections. An equilibrium  $x$ represents, in this context, a set of payments that clear the network in a consistent way. 
A key question is to understand the extent to which a shock hitting the value of the assets of a single node $i$ (perturbation of $c_i$) reflects on the entire network and leads to possible cascade effects. In particular, a default node is defined as one for which the quantity $\sum_jP_{ji}x_j+c_i$ (representing the liquidity of the entity $i$) is below the value of its obligation $w_i$ and the default is called partial or total if, respectively, $\sum_jP_{ji}x_j+c_i>0$ or not. Despite its apparent simplicity, this framework has proved to be very useful for analyzing how losses propagate through the financial system. Previous works including \cite{Liu2010,Acemoglu2015a,Glasserman2015,Ren2016} have analyzed conditions for uniqueness of the clearing payment  equilibrium $x$ and studied its dependence on the exogenous flow vector $c$. In particular, Eisenberg and Noe themselves  \cite{Eisenberg.Noe:2001} find sufficient conditions for uniqueness of clearing payment equilibria $x$ in the special case of non-negative exogenous flow vector $c$ and prove monotonicity and concavity of $x$ as a function of $c$. Glasserman and Young \cite{Glasserman2015} also consider the case of non-negative exogenous flow $c$ and  extend the sufficient conditions for uniqueness of the clearing payment  equilibrium $x$ in \cite{Eisenberg.Noe:2001} to cover the case where the matrix $P$ has spectral radius $\rho(P)<1$. They also estimate the extent to which interconnections increase expected losses and defaults under a wide range of shock distributions, providing bounds on the potential magnitude of network effects on contagion and loss amplification.  \cite{Acemoglu2015a} consider a particular case of the Eisenberg and Noe model where the network is regular and prove that the clearing payment equilibrium is generically unique with respect to values of the exogenous flow $c$ in $\R^n$. Furthermore, they prove rigorous results about the resilience of different network topologies depending on the shock magnitude.  Liu and Statum \cite{Liu2010} use linear programming to provide a sensitivity analysis of Eisenberg and Noe model with respect to certain parameters. Finally,  Ren {it et al.} \cite{Ren2016} explore several sufficient conditions for uniqueness of the clearing payment equilibrium, in particular showing that this holds true in the case where at least one entry of the maximal equilibrium is saturated at its upper bound or at least one entry of its maximal equilibrium is saturated at $0$. 

In quadratic network games, the entries of the vector $x$ represent the actions strategically chosen by $n$ players, each one seeking to maximize a utility function 
$u_{i}\left(x\right)=c_ix_i-{x_{i}^{2}/}{2}+x_{i}\sum_{j} P_{ji} x_{j}$ given by the difference between a linear return and a quadratic cost depending only on her own action plus a bi-linear term coupling her action with those of her neighbors in the network. Here, the entries of the exogenous flow $c$ represent the constant marginal benefits of the individual players from their own actions and coincide with their optimal choices in the absence of network interaction, whereas the nonzero entries of the matrix  $P$ correspond to either strategic complementarities (if they are positive) of strategic substitutes (if they are negative) between neighbor players in the network. In the absence of any constraints on their actions, the players' best responses  are linear functions and Nash equilibria are solutions of the linear system $x=P^{\top}x+c$ whose existence and uniqueness can be characterized in terms of the spectral properties of $P$. In particular, if $P$ has spectral radius $\rho(P)<1$, then there exists a unique Nash equilibrium and, in the case when all externalities are positive, \cite{Ballester.ea:2006} show how its aggregate performance can be evaluated in terms of the sum of the individual players' marginal benefits weighted by their so-called Bonacich network centrality \cite{Bonacich:1987}. 
When upper and lower bounds on the feasible players' actions are considered, the best responses prove to be described as the composition of linear functions with saturation non-linearities and Nash equilibria 
%
coincide with the solutions of the fixed point equation \eqref{eq:nonlinear-system} \cite{CalvoArmengol.ea:2009,Bramoulle.Kranton:2016}. 
In this case, it is known that, while existence is ensured by convexity and compactness of the strategy profile space as argued before, uniqueness is lost in general. In this regard,  \cite{Ballester.ea:2006} claim that 
``multiple equilibria will certainly emerge, which is a plausible outcome in the school setting'', 
while \cite{Bramoulle.Kranton:2016} acknowledge that ``our general knowledge of how unique versus multiple equilibria depend on parameters and the network is still very fragmented.'' 
For symmetric quadratic games of strategic substitutes (i.e., non-positive symmetric $P$),  
Bramoull\'e {\it et al} \cite{Bramoulle.ea:2014} prove uniqueness of Nash equilibria when $P$ has spectral radius $\rho(P)<1$,  building on the fact that in this case the quadratic game is  potential \cite{Monderer.Shapley:1996} with strictly concave potential function.
On the other hand, in the special case when the exogenous flow $c$ is strictly positive, Belhaj {\it et al} \cite{Belhaj.ea:2014} provide sufficient conditions for uniqueness of Nash equilibria for a class of network games with strategic complements (non-negative $P$) that include quadratic games, generalizing a  previous result for fixed points of monotone concave functions \cite{Kennan:2001}.

The present paper develops a systematic study of the network equilibria described by equation \eqref{eq:nonlinear-system} in the general case of networks $(P,w)$ where $P$ is a non-negative square matrix with spectral radius $\rho(P)\le1$ and provides three fundamental contributions:
\begin{enumerate}
	\item[(i)] We characterize a class of non-expansive networks (c.f.~Definition \ref{def:nonexpansive}) including as a special case networks where $P$  is a row-stochastic or sub-stochastic matrix and we prove that, for this class, all network equilibria satisfy an invariance property (Theorem \ref{theorem:partition}) with respect to a certain partition of the node set in  surplus, exposed, and deficit nodes (c.f.~Section \ref{sec:invariance});
	\item[(ii)] We analyze the structure of the set of network equilibria with respect to topological properties of the network. In particular, we show how to effectively construct all network equilibria starting from anyone of them and  prove necessary and sufficient conditions for uniqueness of the network equilibrium in the general case of spectral radius $\rho(P)\le1$ (Theorem \ref{theorem:general}).  This result subsumes and extends the ones available in the previously surveyed literature on  financial networks, as in this context $P$ is always a stochastic or sub-stochastic matrix, hence  with spectral radius $\rho(P)\le1$.  It is worth emphasizing that uniqueness conditions we derive can be easily checked a priori without the need for computing the network equilibrium itself.
\item[(iii)] We show that network equilibria exhibit a jump discontinuity in their dependence on the exogenous flow vector $c$ when this is crossing certain regions of measure $0$ representable as graphs of continuous functions,  where the uniqueness of equilibrium is lost (Theorem \ref{theo:continuity}). 
We provide an analytical description of the discontinuity set and we quantify the size of the largest jump (Corollary \ref{coro:jump}). In the financial network application, this can be interpreted as a jump in the aggregate loss function (c.f.~Section \ref{sec:systemicrisk} and Example \ref{example:loss}).
\end{enumerate}\smallskip
Notice that, in contrast to some of the previously reviewed literature,  we do not make any symmetry or regularity assumptions on the matrix $P$ describing the network (c.f.~\cite{Ballester.ea:2006,CalvoArmengol.ea:2009,Galeotti.ea:2010,Bramoulle.ea:2014,Bramoulle.Kranton:2016,Acemoglu2015,Acemoglu2015a}), nor on the sign of the exogenous flow $c$ (c.f.~\cite{Eisenberg.Noe:2001,Glasserman2015,Belhaj.ea:2014}). This creates several  technical challenges as in particular, we cannot rely on the theory of potential games (which would require $P$ to be symmetric) and we have to deal with possible effective saturations at both the upper and the lower bound (while, e.g.,  assuming non-negative $c$ would have removed the impact of the lower saturation). 

From a methodological viewpoint, it is worth pointing out that non-negativity of the matrix $P$ allows one to interpret the considered network equilibria as the Nash equilibria of a particular class of games with strategic complementarities. This implies that some of the general results in the theory of supermodular games \cite{Topkins:1979,Milgrom.Roberts:1990,Vives:1990,Topkins:1998} can be applied  in order to guarantee, e.g., that the set of network equilibria is a complete lattice, as well as the validity of certain comparative statics \cite{Milgrom.Shannon:1994}, in particular that the minimal and maximal network equilibria are monotone functions of the exogenous flow vector $c$, of the upper saturation vector $w$, and of the matrix $P$ itself (Proposition \ref{proposition:c}). 
However, we depart quite soon from the general theory of supermodular games and develop an approach to the study of such  monotone linear saturated network systems that partly hinges on some of the theory of non-negative matrices \cite{Berman.Plemmons:1994} (cf.~Proposition \ref{proposition:PF}). Key steps in our treatment include the derivation of some {\it ad hoc} technical results exploiting finer spectral and topological properties of the network (Propositions \ref{prop:contraction}, \ref{prop:uniquenss-outconnected}, and \ref{proposition:uniqueness-irreducible}) that then prove instrumental in the proof of our main results (Theorems \ref{theorem:partition}, \ref{theorem:general}, and \ref{theo:continuity}). We notice that our results for non-expansive networks are somewhat reminiscent of the Rural Hospitals Theorem \cite{Roth:1984,Roth:1986} in the matching literature which, under suitable assumptions, shows that the set of stable matchings (hence, the equilibria in that setting) is a distributive lattice and satisfies a fundamental invariance property.

The rest of this paper is organized as follows. The remainder of this Introduction is devoted to a brief explanation of the main notational conventions to be followed throughout the paper.  Section \ref{sec:motivations} presents the two main motivating applications for the model considered, i.e., financial networks and network games with linear saturated best replies. Section \ref{sec:structure} establishes a number of preliminary results on the structure of the equilibria. Uniqueness results as well a general expression describing all solutions in non-uniqueness cases is presented in Section \ref{uniqueness}. Section \ref{discontinuous} is devoted to the analysis of jump discontinuities in the equilibrium with respect to the variation of the exogenous flow vector with a particular focus on financial networks. The paper ends with Section \ref{conclusions} dedicated to draw some conclusions and open problems.
\medskip

\textbf{Notation}
We briefly explain the notation to be used throughout this paper. 
Vectors are denoted with lower case, matrices with upper case, and sets with calligraphic letters. A subscript associated to vectors, for instance $v_{\mc A}$, represents the sub-vector that is the restriction of a vector $v$ in $\R^n$ on the set of indices $\mc A \subseteq \{1,2, \dots, n\}$. The same notation is used for matrices: $P_{\mc A \mc B}$ represents the sub-matrix of $P$ obtained by considering rows and columns associated with the indices contained in sets $\mc A$ and $\mc B$, respectively. We indicate with $\1$ the all-$1$ vector, regardless of its dimension. 

Throughout the paper, the natural entry-wise partial order is considered on $\R^n$,  so that, the inequality $ x\leq y$ for two vectors $x$ and $y$ in $\R^n$ is to be understood as  $x_i\leq y_i$ for every $i=1,2,\ldots,n$, whereas $x\lneq y$ means that $x\le y$ with strict inequality in at least one entry. Analogously, the absolute value of a vector $v$ in $\C^n$ is the vector $|v|$ in $\R_+^n$ with entries $(|v|)_i=|v_i|$ for $i=1,\ldots,n$. A norm $\|\,\cdot\,\|$ on $\C^n$ is referred to as monotone if $\|v\|\le\|w\|$ whenever $|v|\le|w|$. 

%

The \emph{spectral radius} of a square matrix $P$ in $\R^{n\times n}$ is denoted by $\rho \left(P\right)$. 
A non-negative square matrix $P$ in $\mathbb{R}_{+}^{n \times n}$ is referred to as \emph{(row) sub-stochastic} if the sum of the entries in each row never exceeds $1$, i.e., if $P\1 \le  \1$. Notice that in the literature it is often assumed that sub-stochastic matrices have the additional property that for at least one row there is strict inequality. Here we prefer not to follow this convention and in this way our class of sub-stochastic matrices contains also the \emph{stochastic} matrices that are those for which $P\1=\1$.   A non-negative square matrix $P$  is irreducible if for every $i$ and $j$, there exists $l\ge1$ such that $(P^l)_{ij}>0$. 

A directed graph is the pair of a finite node set $\mc V$ and of a set $\mc E\subseteq\mc V\times\mc V$ of links, whereby a link $(i,j)$ is meant as directed from its tail node $i$ to its head node $j$. To any square matrix $P$ in $\R^{n\times n}$, we associate a directed graph $\mc G_P=(\mc V,\mc E)$ with node set $\mc V=\{1,2,\ldots,n\}$, and link set $\mc E=\{(i,j)\in\mc V\times\mc V:\,P_{ij}\ne0\}$.  




\section{Applications}\label{sec:motivations}
In this section, we describe two main motivating applications. We start in Section \ref{sec:financial} by presenting a model of financial networks generalizing the one first considered in \cite{Eisenberg.Noe:2001}. We then provide an interpretation of network equilibria as Nash equilibria for a class of network games with monotone linear saturated best responses, as explained in Section \ref{sec:game}.
Notice that the considered notion of equilibrium in saturated networks and the results derived in the following sections  may find application in other contexts, such as, e.g., in some dynamical network flow models \cite{Massai.ea:IFAC2020}. 


\subsection{Payment equilibria in financial networks}\label{sec:financial}
We consider a set $\mc V=\{1,\dots , n\}$ of financial entities (e.g., banks, broke dealers,...) interconnected by internal and external obligations that are specified by a non-negative matrix $W$ in $\R_+^{n\times n}$ and three non-negative vectors $a$, $b$, and $u$ in $\R^n_+$ whose entries have the following interpretation:
\begin{itemize}
	\item $W_{ij} \ge 0$ is the
	liability of node $i$ to node $j$;
	\item $a_i$ is the total value of assets and credits of $i$
	from external entities;
	\item $b_i$ is the total liability
	of node $i$ to external non-financial entities;
	\item $u_i$ is the total liability
	of node $i$ to external financial entities. 
\end{itemize}
The quantity $v_i=\sum_j W_{ji}-\sum_j W_{ij}+a_i-b_i-u_i$ is the net worth of node $i$. If the condition $v_i\geq 0$ is verified for every $i$ in  $\mc V$, it means that each node is fully liable and in principle capable to pay back all its liabilities to the nodes in the network as well the external ones. In case when instead some nodes do not satisfy the condition $v_i\geq 0$, namely they are not fully liable, it is necessary to determine a consistent set of payments among the various nodes. 

Put $w_i=\sum_j W_{ij}+u_i$ and 
$$P_{ij}=\left\{\begin{array}{ll}W_{ij}/w_i\quad \quad  &\text{if}\, \: \:  w_i>0\\[7pt] 0\quad &\text{otherwise}\end{array}\right.$$
We define by $X_{ij}$ the payment from node $i$ to node $j$ and by $X_{i{\rm o}}$ the payment from node $i$ to external financial entities. Assuming that liabilities to non-financial entities have a higher seniority and that all other payments (including those to external financial entities) should be proportional to the corresponding liabilities, a consistent set of payments among the nodes has to satisfy the relations
\begin{equation}\label{payments}
\begin{array}{rcl}X_{ij} &=&\ds  \min \left\{  P_{ij} \max \left\{ \sum_{k} X_{ki}+a_i-b_i,0 \right\}, W_{ij}\right\}\\[15pt]
X_{i{\rm o}} &=&  \min \left\{  \ds\frac{u_i}{w_i} \max \left\{ \sum_{k} X_{ki}+a_i-b_i,0 \right\}, u_i\right\}\end{array}
\end{equation}
%
Let $x_i=\sum_{j }X_{ij}+X_{i{\rm o}}$ be the total payment of node $i$ to the financial entities. Summing the relations in \eqref{payments} and using the fact that $W_{ij}=w_iP_{ij}$, we obtain
\begin{equation}\label{eq.a0}
x_i=\min \left\{  \max \left\{ \sum_{k } X_{ki}+a_i-b_i,0 \right\}, w_i\right\}
\end{equation}
so that, $X_{ij}=x_iP_{ij}$. Relation \eqref{eq.a0} can thus be rewritten as
\begin{equation}\label{eq.a1}
x_i=\min \left\{  \max \left\{ \sum_{k } x_kP_{ki}+a_i-b_i,0 \right\}, w_i\right\}
\end{equation}
This set of relations is equivalent to \eqref{payments}. Indeed, if the vector $x$ solves \eqref{eq.a1}, then $X_{ij}=x_iP_{ij}$ solves \eqref{payments}.
This coincides with \eqref{eq1} with exogenous flow $c=a-b$. It is worth noticing that, in the financial jargon, vectors $x$ are called \emph{clearing vectors}. 

Notice that the matrix $P$ is sub-stochastic in its strict sense (i.e., at least one row does not sum to $1$) when either there exist nodes with a positive liability towards external financial entities, or nodes with no financial liabilities.

In this financial setting, it is often considered the case when we start from a fully liable configuration, that is $v_i\geq 0$ for all $i$, leading to a solution $x$ of \eqref{eq.a1} such that $x_i\geq w_i$ for all $i$. We then assume that the outside assets suffer a shock $\epsilon\in \mathbb{R}^n_+$ so that their values reduce to $a-\epsilon$ possibly making some of the $v_i$'s negative. The study of the number of nodes in default $x_i< w_i$ as a function of the shock $\epsilon$ is one of the key issues.

{\subsection{Network games with monotone linear saturated best responses} \label{sec:game}
We consider games with player set $\mc V=\{1,\dots , n\}$, whereby each player $i$  in $\mc V$ chooses an action $x_i$ from the compact interval $\mc A_i=[0, w_i]$, where $w_i>0$. We gather all actions in a vector $x$ to be referred to as the strategy profile. Following a standard notational convention in game theory, we indicate by $x_{-i}$ in $\prod_{j\ne i}\mc A_j$ the strategy profile of all players other than player $i$. 

First, we consider the case of quadratic utility functions 
\begin{equation}\label{utilities:quadratic}u_{i}\left(x\right)=u_i(x_i, x_{-i})=c_ix_i-\frac{x_{i}^{2}}{2}+x_{i}\sum_{j} P_{ji} x_{j}\,,\end{equation}
for every player $i$ in $\mc V$ and strategy profile $x$. In \eqref{utilities:quadratic}, $c_i$ denotes the marginal benefit of individual $i$ from its own action, while $P$ is a non-negative matrix describing the strategic interactions among the various players. Notice that, absent network effects, i.e., in the special case $P=0$, $c_i$ is the optimal action of player $i$.  

Such games are known in the literature as constrained quadratic network games. 
Notice that the quadratic utility function $u_i$ in \eqref{utilities:quadratic} implies that the best response of a player $i$ in $\mc V$ is always unique and given by
\begin{equation}\label{br-supermodular}
	B_i(x_{-i})= \min \left\{ \max \left\{ \sum_{j=1}^n x_jP_{ji}+c_i, 0 \right\} , w_i \right\}\,. 
	\end{equation}
It follows that Nash equilibria for such constrained quadratic network games are exactly the solutions of the fixed point equation \eqref{eq:nonlinear-system}.

In this paper, we focus on the special case where the coefficients $P_{ji}$ are all non-negative. In this way, we are considering games of pure strategic complements: for every player $i$, the higher the value of $x_{-i}$, the higher the rate of variation of the utility $u_i(x_i, x_{-i})$ of player $i$  with respect to its own action $x_i$. Mathematically, games like these, where actions belong to compact spaces and utilities $u_i$ are twice differentiable functions with non-negative cross derivatives 
$$\frac{\partial^2u_i}{\partial x_i\partial x_j}=P_{ji}\geq 0$$
for every $i$ and $j$ with $j\neq i$, are called \emph{supermodular}.  A more general definition of supermodular games can be found in \cite{Milgrom.Roberts:1990}; the one proposed here will be sufficient for our purposes.
It is known that supermodular games always admit a complete lattice of Nash equilibria and in our case they coincide with the solutions of \eqref{eq:nonlinear-system}. This fact will be exploited in the Section \ref{sec:lattice}. 

In fact, our analysis applies to the broader class of network games with linear saturated best response as in \eqref{br-supermodular}. This includes, e.g., games with player set $\mc V$, action space $\mc A_i=[0,w_i]$,  for every player $i$ in $\mc V$, and utility functions in the form 
\be\label{utility:general1}u_i(x)=\varphi_i\left(x_i-c_i+\sum_{j\ne i}P_{ji}x_j\right)\,,\ee
for a continuous function $\varphi_i:\R\to\R$  that is increasing on $(-\infty,0]$ and  decreasing in $[0,+\infty)$  \cite{Bramoulle.Kranton:2016}. Notice that \eqref{utilities:quadratic} is a special case of \eqref{utility:general1} with $\varphi_i(y)=-y^2/2$.

\section{Structural properties of network equilibria}\label{sec:structure}



While existence of network equilibria is guaranteed for every network $(P,w)$ and exogenous flow $c$, as discussed in Section \ref{sec:intro}, their uniqueness or multiplicity and more generally the structure of the network equilibrium set $\mc X$ remain more delicate issues, as also illustrated in the following simple example. \medskip
\begin{figure}
\begin{center}
\includegraphics[height=4cm]{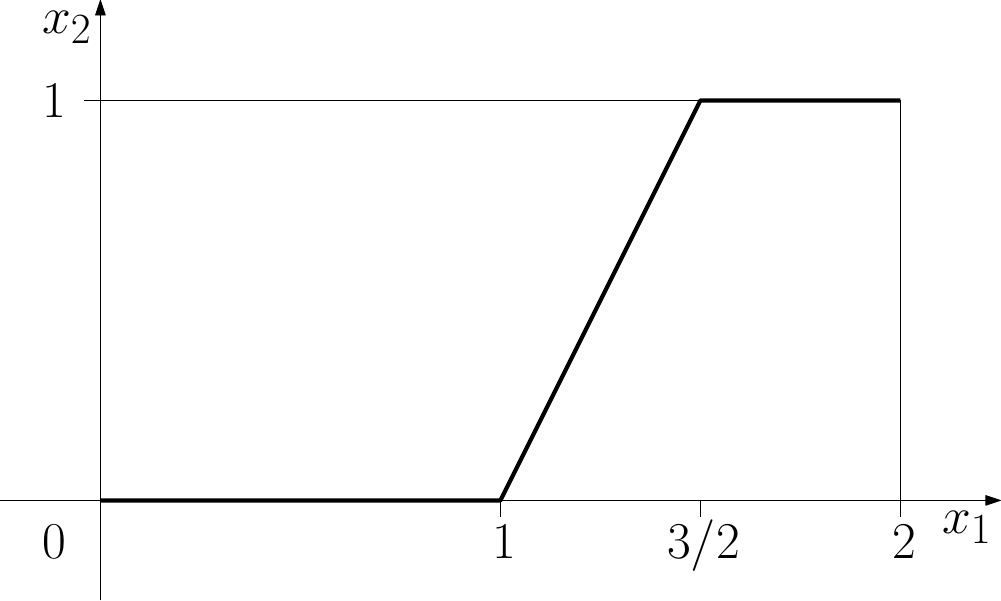}
\caption{Set of network equilibria for the network in Example \ref{example:non-partition}.\label{fig:equilibrium-set}}
\end{center}
\end{figure}

\begin{example}\label{example:non-partition}
Consider a network $(P,w)$ with $n=2$, 
$$P=\left[\ba{cc}1&1\\0&\frac 12\ea\right]\,,\qquad w=\left[\ba{c}2\\1\ea\right]\,,$$
and the exogenous flow
$$c=\left[\ba{c}0\\-1\ea\right]\,.$$
In this case, the fixed-point equation \eqref{eq:nonlinear-system} reads 
\be\label{eq:example}x_1=S_0^2(x_1)\,,\qquad x_2=S_0^1(x_1+x_2/2-1)\,,\ee
and the set of network equilibria is then 
\be\label{eq:equilibrium-set-example}\mc X=\left\{\left(t,S_0^{1}(2t-2)\right):\, 0\le t\le2\right\}\,,\ee
as displayed in Figure \ref{fig:equilibrium-set}. 
\end{example}\medskip

In the rest of this section, 
we study structural properties of the set of network equilibria $\mc X$ for a network $(P,w)$ with exogenous flow $c$, i.e., for the set of solutions of the fixed-point equation \eqref{eq:nonlinear-system}.  {\color{black}Specifically, the contribution of this section is threefold. First, we exploit the fact that the network equilibrium set $\mc X$ can be interpreted as the set of Nash equilibria of the $n$-player supermodular game with utilities as in \eqref{utilities:quadratic} and we establish a number of results concerning the lattice structure of $\mc X$ and its monotone dependence on the exogenous flow  $c$. Second, we review some classical results on the spectral theory of non-negative matrices and derive some additional properties of the set of network equilibria $\mc X$ for a special class of non-expansive networks. Third, we introduce a fundamental partition of the node set into three subsets and prove that such partition is invariant with respect to the entire set of network equilibria for non-expansive networks.  We wish to remark that, while the results concerning the lattice structure hold true in general for every network $(P,w)$, the rest of the results are instead deeply connected to the finer spectral assumptions on the matrix $P$ (c.f.~Definition \ref{def:nonexpansive}) and do not hold true for general networks. In particular, such results involve properties of the network equilibrium set that will play a crucial role in the following sections. 


%
%
\subsection{Lattice properties of the set of network equilibria}\label{sec:lattice}
For a network  $(P,w)$ and an exogenous flow $c$, consider the following recursion on the complete lattice $\mc L_0^w$: 
\be\label{eq:iterative-algorithm}x(t+1)=S_0^w(P^{\top}x(t)+c)\,,\qquad t\ge0\,.\ee
Notice that equation \eqref{eq:iterative-algorithm} can be interpreted as the update rule of a synchronous best response dynamics for the supermodular game with utilities as in \eqref{utilities:quadratic}. 
The following proposition gathers a number of results on the network equilibria set $\mc X$ that follow from \cite{Milgrom.Roberts:1990} as a direct consequence of such game-theoretic interpretation.
\begin{proposition}\label{proposition:c} Consider a network $(P,w)$ and an exogenous flow $c$ and let $\mc X$ be the corresponding set of network equilibria. Let $x(t)$, for $t=0,1,\ldots,$ be the sequence generated by the recursion \eqref{eq:iterative-algorithm} with initial condition $x(0)=x_0$ in $\mc L_0^w$. Then:
	\begin{enumerate}
		\item[(i)] $\mc X$ is a complete lattice in $\R^n$. In particular,  there exist a minimal network equilibrium $\ul x$ and a maximal network equilibrium $\ov x$ in $\mc X$; 
		\item[(ii)] if $x_0=0$, then $x(t)$ is non-decreasing and $\lim x(t)=\ul x$ as $t$ grows large; 
		\item[(iii)] if $x_0=w$, then $x(t)$ is non-increasing and  $\lim x(t)=\ov x$ as $t$ grows large;
		\item[(iv)] both $\ul x$ and $\ov x$ are monotone non-decreasing functions of the exogenous flow $c$ in $\R^n$, of the matrix $P$ in $\R_+^n$, and of the upper saturation vector $w$ in $\R_+^n$.
	\end{enumerate}
\end{proposition}
}

\begin{remark}\label{remark:algorithm}
Observe that the recursion \eqref{eq:iterative-algorithm} can be implemented as a distributed iterative algorithm, whereby at every time $t=0,1,\ldots$, each node $i$ in $\mc V$ updates in parallel its state $x_i(t)$ to $$x_i(t+1)=S_0^{w_i}\left(\sum\nolimits_j P_{ji}x_i(t)+c_i\right)\,.$$ Notice that such update only requires each node $i$ to observe the current states $x_j(t)$ of its in-neighbors $\{j\in\mc V:\,P_{ji}>0\}$ and the total complexity of each iteration of \eqref{eq:iterative-algorithm}  is of the order of the number of links in the network, i.e., the number of non-zero entries of $P$. 

 We now make some more refined considerations on the convergence time.
Consider the recursion \eqref{eq:iterative-algorithm} with the initial condition $x(0)=0$ and let $t^-_i=\inf\{t\ge0:\,x_i(t)=w_i\}$ for every $i=1,\ldots,n$. By Proposition \ref{proposition:c} (ii), whenever $t_i^-<+\infty$ we have $x_i(t)=w_i$ for every $t\ge t_i^-$.
Analogously, by considering the recursion \eqref{eq:iterative-algorithm} this time with the initial condition $x(0)=w$ and letting $t^+_i=\inf\{t\ge0:\,x_i(t)=0\}$ for $i=1,\ldots,n$, Proposition \ref{proposition:c} (iii) guarantees that, whenever $t_i^+<+\infty$ we have $x_i(t)=0$ for every $t\ge t_i^+$. Observe that, since $\ov x\ge\ul x$, we necessarily have 
that at most one between (and possibly neither of) $t_i^-$ and $t_i^+$ is finite.  Let  $t_i=\min\{t_i^-,t_i^+\}$ for all $i=1,\ldots, n$. Then, when $t^*=\max\{t_i:\,1\le i\le n\}<+\infty$, we have a unique network equilibrium $x^*=\ov x=\ul x$ with every entry saturated from either below or above and  convergence in finite time $t^*$ is guaranteed to $x^*$ from every initial condition $x(0)$ in $\mc L_0^w$.  In contrast, when $t_i=+\infty$ for some $i$ convergence typically occurs in infinite time, see Remark \ref{remark:finitetime} for further considerations in this case. 
\end{remark}\medskip

\subsection{Spectral properties and non-expansive networks} \label{sec:non-expansive}
In this subsection, we derive a number of notions and results on non-negative matrices and introduce the notion of non-expansive network that will play a key role in the rest of the paper. 
We start with the following proposition gathering known results that can be found, e.g., in the monograph \cite{Berman.Plemmons:1994}. 

\begin{proposition}\label{proposition:PF} Let $P$ in $\R_+^{n\times n}$ be a non-negative square matrix. Then:
\begin{enumerate} 
\item[(i)] the spectral radius $\rho(P)$ is an eigenvalue of $P$ and there exist vectors $p$ and $\pi$ in $\R_+^n\setminus\{0\}$ such that $$Pp=\rho(P)p\,,\qquad \pi^{\top}P=\rho(P)\pi^{\top}\,.$$ Such vectors are called, respectively, a right and a left dominant eigenvector of $P$;
\item[(ii)] if $Q$ is a principal square sub-matrix of $P$, then $\rho(Q)\leq \rho(P)$.\end{enumerate}
Moreover, if $P$ is irreducible, then 
\begin{enumerate}
\item[(iii)] the dominant eigenvectors $p$ and $\pi$ are unique up to normalization and have all positive entries;
\item[(iv)]  for every vector $c$ in $\R^n$ such that $p^{\top}c=0$, the equation $x=\rho(P)P^{\top}x+c$ admits infinite solutions $x$ in $\R^n$;
\item[(v)]  if $Q$ is a principal proper square sub-matrix of $P$, then $\rho(Q)< \rho(P)$.
\end{enumerate}
\end{proposition}

For a non-negative square matrix $P$ in $\R_+^{n\times n}$, we shall consider the connected components $\mc V_1,\dots ,\mc V_s$ of the associated digraph $\mc G_P$ and refer to them as the \emph{classes} of $P$.
Upon a possible permutation of the indices $i=1,\ldots,n$, we can always  assume that the matrix $P$ admits the following block triangular structure 
\begin{equation}\label{canonical} P=\left[\begin{matrix} P^{(11)} &P^{(12)} &\cdots &P^{(1s)}\\ 0 & P^{(22)}& \cdots & P^{(2s)}\\ 0&0& \ddots & \vdots\\ 0&0& \cdots &P^{(ss)}\end{matrix}\right]\,,\end{equation}
where, for $i,j=1,\ldots,l$, $P^{(i,j)}$ in $\R_+^{\mc V_i\times\mc V_j}$ is the sub-matrix of $P$ obtained by keeping only rows with index in $\mc V_i$ and columns with index in $\mc V_j$. Notice that this is equivalent to saying that the diagonal blocks $P^{(ii)}$ are irreducible and that in $\mc G_P$ there is no link from a node in a class $\mc V_l$ to any node in a class $\mc V_i$ with $i<l$. 
It then follows from Proposition \ref{proposition:PF} (ii) that $\rho(P^{(ii)})\leq\rho(P)$.  
A class $\mc V_i$, for $1\le i\le s$, will then be referred to  (c.f.~\cite{Berman.Plemmons:1994}) as:  
\begin{itemize}
\item\emph{basic} if $\rho(P^{(ii)})=\rho(P)$; 
\item\emph{final} if $P^{(ih)}=0$ for every $h\ne i$.
\end{itemize}

We can state the following result. 
\begin{proposition}\label{prop:contraction} Let $P$ in $\R_+^{n\times n}$  be a non-negative square matrix. Then, there exists a positive vector $v$ in $\R^n_+$ such that \be\label{Pv<=v}Pv\leq v\,,\ee
if and only if $\rho(P)< 1$ or $\rho(P)=1$  and every basic class of $P$ is final.
\end{proposition}
\begin{proof} See Appendix \ref{proof:contraction}.\end{proof}\medskip

Observe that to every positive vector $v$ in $\R^n_+$ we may associate the weighted $l_1$-norm
\be\label{eq:v-norm}\|x\|=\sum_{i=1}^nv_i|x_i|\,,\qquad x\in\C^n\,.\ee
Clearly, the above is an absolute norm, hence a monotone  norm \cite{Johnson.Nylen:1991}. 
Condition \eqref{Pv<=v} implies that 
\be\label{eq:non-expansive}\|P^{\top}x\|=v^{\top}P^{\top}|x|\le v^{\top}|x|=\|x\|\,,\qquad \forall x\in\C^n\,.\ee
We introduce the following definition. 
\medskip
\begin{definition}\label{def:nonexpansive}
A network $(P,w)$ is \emph{non-expansive}  if either \begin{enumerate}
\item[(i)] $\rho(P)< 1$; or 
\item[(ii)] $\rho(P)=1$  and every basic class of $P$ is final.
\end{enumerate} 
 \end{definition}

\medskip

\begin{remark} \label{remark:non-expansive1}
Consider a \emph{non-expansive} network $(P,w)$ and let $\|\,\cdot\,\|$ be the monotone vector norm defined by \eqref{eq:v-norm} for a positive vector $v$ satisfying \eqref{Pv<=v}. Then, for arbitrary vectors $x,\tilde x,c,\tilde c$ in $\R^n$, we have 
\be\label{eq:remark-non-expansive}\ba{rcl}\|S_{0}^w(P^{\top}x+c)-S_{0}^w(P^{\top}\tilde x+\tilde c)\|
&=& \ds\sum_{i=1}^nv_i|S_{0}^{w_i}((P^{\top}x)_i+c_i)-S_{0}^{w_i}((P^{\top}\tilde x)_i+\tilde c_i)|\\
&\le& \ds\sum_{i=1}^nv_i|(P^{\top}x)_i+c_i-(P^{\top}\tilde x)_i-\tilde c_i|\\[12pt]
&\le& \ds\sum_{i=1}^nv_i|(P^{\top}(x-\tilde x))_i|+\sum_{i=1}^nv_i|c_i-\tilde c_i|\\[12pt]
&=& \|P^{\top}(x-\tilde x)\|+\|c-\tilde c\|\\[7pt]
&\le&\|x-\tilde x\|+\|c-\tilde c\|\,,
\ea\ee
where the first inequality above follows from monotonicity of the weighted $l_1$-norm $\|\,\cdot\,\|$ and the last one from \eqref{eq:non-expansive}. This property justifies the terminology introduced in Definition \ref{def:nonexpansive}. 
\end{remark}\medskip

\begin{remark} \label{remark:non-expansive2} 
A special class of non-expansive networks is provided by those networks $(P,w)$ such that the matrix $P$ is \emph{ (row) sub-stochastic}, that is a matrix $P \in \mathbb{R}_{+}^{n \times n}$ where the sum of the elements in each row never exceeds $1$, i.e., $P\1 \le  \1$. 
Indeed, for a sub-stochastic matrix $P$, it can be easily checked that $\rho(P)\le 1$ and that if $\rho(P)=1$ then every basic class is necessarily final. 
Notice that in the literature it is often assumed that sub-stochastic matrices have the additional property that for at least one row there is strict inequality. Here we prefer not to follow this convention and in this way our class of sub-stochastic matrices contains also the \emph{stochastic} matrices that are those for which $P\1=\1$. 
\end{remark}\medskip

\begin{remark} \label{remark:non-expansive3}
It is worth pointing out that existence of a (not necessarily monotone) vector norm $\|\,\cdot\,\|$ on $\C^n$ such that \eqref{eq:non-expansive} holds true can be guaranteed under slightly weaker assumptions than those in Definition \ref{def:nonexpansive}. Specifically \cite{Lim:2009} shows that this is equivalent to that either $\rho(P)<1$ or $\rho(P)=1$ and the geometric multiplicity of every eigenvalue $\lambda$ of $P$ with $|\lambda|=1$ is equal to its algebraic multiplicity. In fact, notice that, when $\rho(P)=1$, that every basic class of $P$ is final implies that the geometric multiplicity of every eigenvalue $\lambda$ of $P$ with $|\lambda|=1$ is equal to its algebraic multiplicity, but not vice versa. For a counterexample, take $P$ as in Example \ref{example:non-partition}:  there $P$ has unitary spectral radius and $\lambda=\rho(P)=1$ is a simple eigenvalue, with algebraic and geometric multiplicities both equal to $1$, however, there are two classes, $\mc V_1=\{1\}$  and $\mc V_2=\{2\}$, the first of which is basic but not final. 

In fact, such stricter condition (ii) in Proposition \ref{prop:contraction} in the case when $\rho(P)=1$ ensures not only existence of a vector norm $\|\,\cdot\,\|$ on $\C^n$ such that \eqref{eq:non-expansive} holds true, but also that such a vector norm can be chosen as a weighted $l_1$-norm \eqref{eq:v-norm}. It is exactly the monotonicity of such a norm that allows one to show that non-expansiveness is preserved when composing the affine map $P^{\top}x+c$ with the nonlinear saturation $S_0^w(\,\cdot\,)$, as in \eqref{eq:remark-non-expansive}. 
\end{remark}


\subsection{Invariance property of network equilibria}\label{sec:invariance}
In this subsection, we show that the set of network equilibria $\mc X$ of every non-expansive network presents a relevant invariant property that will play a key role in the uniqueness results presented in the next section.

Consider an arbitrary network $(P,w)$ with exogenous flow $c$. For a network equilibrium $x$ in $\mc X$, we can always introduce the following partition of the node set $\mc V=\{1,2,\ldots,n\}$: 
\be\label{eq:node-partition}\mc V=\mc V^x_- \cup \mc V^x_+ \cup \mc V^x_0\,,\ee
where
\begin{itemize}
	\item $\mc V^x_+=\l\{i\in\mc V:\,c_i+\sum_{k \ne i} P_{ki}x_k>w_i\r\}$ is the set of \emph{surplus} nodes; 
	\item $\mc V^x_0=\l\{i\in\mc V:\,0 \le c_i+\sum_{k \ne i} P_{ki}x_k \le w_i\r\}$ is the set of \emph{exposed} nodes; 
	\item $\mc V^x_-=\l\{i\in\mc V:\,c_i+\sum_{k \ne i} P_{ki}x_k<0\r\}$ is the set of  \emph{deficit} nodes. 
\end{itemize} 
Observe that, by the way these sets have been defined, it directly follows that
\be\label{eq:deficit-surplus}
\ba{ll}x_i=0\,,\quad & \forall\, i\in\mc V^x_-\,,\\[7pt]
x_i=w_i\,,\quad\; & \forall\, i\in\mc V^x_+\,, \\[7pt]
x_i=c_i+\sum_{j \ne i} P_{ji}x_j\,, \quad\; &\forall\, i\in\mc V^x_0\,.\end{array}\ee

We now show that, if the network $(P,w)$ is non-expansive, then partition \eqref{eq:node-partition} is invariant with respect to the chosen network equilibrium. 
This is stated in the following, which is the key result of this section and will be instrumental to all our future derivations.

\begin{theorem}\label{theorem:partition} 
For a non-expansive network $(P,w)$ and any exogenous flow $c$ in $\R^n$, the partition \eqref{eq:node-partition} is invariant over all equilibria $x$ in $\mathcal X$.
	\end{theorem}
\proof
We shall consider the maximal network equilibrium $\ov{x}$ and any another network equilibrium $x$ in $\mc X$ and show that they share the same node partition \eqref{eq:node-partition}. To begin with, notice that, since $\ov x\ge x$, we have $\mc V^{\ov{x}}_+ \supseteq \mc V^x_+ $ and $ \mc V^{\ov{x}}_- \subseteq \mc V^x_-$. Let us split nodes in five different classes, $\mc C_1, \mc C_2, \mc C_3, \mc C_4, \mc C_5$, corresponding to the possible cases in which the entries of the network equilibria $\ov{x}$ and $x$ can differ and are precisely defined as follows:
\begin{itemize}
	\item $\mc C_1 = \mc V^x_+$ is the set of nodes that are surplus for both equilibria; 
	\item  $\mc C_2=\mc V^{\ov{x}}_+ \setminus \mc V^x_+$ is the set of nodes that are surplus for $\ov{x}$ but not for $x$; 
	\item $\mc C_3 = \mc V^x_0 \cap \mc V^{\ov{x}}_0$ is the set of nodes that are exposed for both equilibria; 
	\item $\mc C_4 = \mc V^{\ov{x}}_0 \setminus \mc V^x_0$ is the set of nodes that are exposed for $\ov{x}$ and deficit for $x$; 	
	\item $\mc C_5 = \mc V^{\ov{x}}_-$ is the set of nodes that are deficit for both equilibria. 
	\end{itemize}	
We shall  write any vector $y$ in $\R^n$ in a block form $y=(y^{(1)}, y^{(2)},y^{(3)},y^{(4)}, y^{(5)})$ and for simplicity of notation indicate $Q^{(ij)}:=\left(P^{\top}\right)_{\mc C_i\mc C_j}$ for $i,j=1,\ldots,5$.
Notice that 
$\ov x^{(1)}=x^{(1)}=w^{(1)}$, $\ov{x}^{(5)}=x^{(5)}=0$, and  
	\begin{equation}\label{rel1}w^{(2)}=\ov{x}^{(2)}<\sum_{k =1 }^4 Q^{(2k)}\ov{x}^{(k)}+c^{(2)},\qquad x^{(2)} \ge \sum_{k =1 }^4 Q^{(2k)}x^{(k)}+c^{(2)}\,,\end{equation}
	\begin{equation}\label{rel2}\ov{x}^{(3)}=\sum_{k =1 }^4 Q^{(3k)}\ov{x}^{(k)}+c^{(3)},\qquad x^{(3)}=\sum_{k =1 }^4 Q^{(3k)}x^{(k)}+c^{(3)}\,,\end{equation}
	\begin{equation}\label{rel3}\ov{x}^{(4)}=\sum_{k =1 }^4 Q^{(4k)}\ov{x}^{(k)}+c^{(4)},\qquad 0=x^{(4)}>\sum_{k =1 }^4 Q^{(4k)}x^{(k)}+c^{(4)}\,.\end{equation}

Put $z=\ov{x}-x \ge 0$ and notice that, for classes $\mc C_1$ and $\mc C_5$ we have that $z^{(1)}=z^{(5)} =0$. For the remaining blocks, using \eqref{rel1}, \eqref{rel2}, and \eqref{rel3}, we obtain
	\begin{equation}
\label{key-1}
	z^{(2)}  
	< 
	\sum_{k =2 }^4 Q^{(2k)}z^{(k)}\,,\qquad 
	z^{(3)}  
	= 
	\sum_{k =2 }^4 Q^{(3k)}z^{(k)}\,,\qquad 
	z^{(4)} 
	< 
	\sum_{k =2 }^4 Q^{(4k)}z^{(k)}\,.
	\end{equation}
Now, assume by contradiction that $\mc C_2\cup \mc C_4\neq\emptyset$, so that  the above would imply that 
\be\label{z<=P'z}z\lneq P^{\top}z\,.\ee
Since the network is non-expansive, by Proposition \ref{prop:contraction} there exists a positive vector $v$ such that \eqref{Pv<=v} holds true. Together with \eqref{z<=P'z}, this would imply that 
$$v^{\top}z<v^{\top}P^{\top}z\le v^{\top}z\,,$$
%
thus leading to a contradiction. This implies that necessarily $\mc C_2=\mc C_4=\emptyset$, so that $z=0$, thus showing invariance of the node partition \eqref{eq:node-partition} with respect to the network equilibria $x$ in $\mc X$.
\qed


We gather some immediate consequences of Theorem \ref{theorem:partition} in the following result. 

\begin{corollary}\label{corollary:partition}
Let $(P,w)$ be a non-expansive network. Then, for every exogenous flow $c$, there exists a partition of the node set 
\be\label{eq:partition-unique}\mc V=\mc V_+\cup \mc V_0\cup \mc V_-\,,\ee
such that, indicated with $z=(z^{(+)}, z^{(0)}, z^{(-)})$ the corresponding block decomposition of a vector $z$ in $\R^n$ and with $P^{(\alpha\beta)}=P_{|\mc V_\alpha\times\mc V_\beta}$ for $\alpha,\beta=-,0,+$, 
\begin{enumerate}
\item[(i)] for every network equilibrium $x$ in $\mc X$ 
\be\label{not0uniqueness2} x^{(-)}=0\,,\qquad x^{(0)}=P^{(00)\top}x^{(0)}+P^{(+0)\top}x^{(+)} +c^{(0)}\,,\qquad x^{(+)}=w^{(+)}\,;\ee
\item[(ii)] for every two network equilibria $x$ and $y$ in $\mc X$, 
\be\label{not0uniqueness} x^{(-)} =y^{(-)}\,,\quad x^{(+)} =y^{(+)}\,.\ee
\end{enumerate}
\end{corollary}

Corollary \ref{corollary:partition} implies that uniqueness can always be tested by simply looking at those entries of the network equilibria that belong to $\mc V_0$ and that such entries solve a linear system of equations. However, the outstanding difficulty in the analysis of the equilibrium set $\mc X$ stems from the fact that the partition \eqref{eq:partition-unique}  is not known a priori, a problem that will be dealt with in the next section. 

\medskip

\begin{remark}
The necessity of the additional assumption that every basic class of $P$ is final for networks $(P,w)$ where $P$ is non-stochastic and $\rho(P)=1$ is illustrated by Example \ref{example:non-partition}. In the network considered there, $P$ has two classes: $\{1\}$ that is basic but not final and $\{2\}$ that is final but not basic. In fact, it is easily seen from \eqref{eq:example} and \eqref{eq:equilibrium-set-example} that, while node $1$ is always exposed for every network equilibrium $x$ in $\mc X$, node $2$ is: 
\begin{itemize} 
\item a deficit node for every network equilibrium $x$ in  $\mc X_-=\{(t,0):\, 0\le t<1\}$; 
\item an exposed node for every network equilibrium $x$ in  $\mc X_0=\{(t,2t-2):\, 1\le t\le 3/2\}$; 
\item  a surplus node for every network equilibrium $x$ in  $\mc X_+=\{(t,1):\, 3/2< t\le 2\}$.
\end{itemize} 
Therefore, partition \eqref{eq:node-partition} is clearly equilibrium-dependent in this case. As already pointed out in Remark \ref{remark:non-expansive3} in this case the matrix $P$ has unitary spectral radius and its eigenvalue $\lambda=\rho(P)=1$ has algebraic and geometric multiplicities both equal to $1$. This shows that, when $\rho(P)=1$, the weaker condition that $\lambda=1$ has algebraic multiplicity equal to its geometric multiplicity  is not sufficient for the conclusions of Theorem \ref{theorem:partition} and Corollary \ref{corollary:partition} to hold true and the stricter assumption that every basic class be final is required. 
\end{remark}\medskip

\begin{remark}\label{remark:finitetime}
For a non-expansive network, consider once again the recursion \eqref{eq:iterative-algorithm} and, for $i=1,\ldots,n$, let $t_i$ be defined as in Remark \ref{remark:algorithm}. Assume that $t_i=+\infty$ for some $i$ and let $t^*=\max\left(\{0\}\cup\{t_i:t_i<+\infty\}\right)$. Then, by combining the considerations in Remark \ref{remark:algorithm} with Theorem \ref{theorem:partition} we get that  the recursion \eqref{eq:iterative-algorithm} started in $x(0)=0$ and $x(0)=w$ respectively determines partition \eqref{eq:partition-unique} by time $t^*$. Indeed, the surplus, deficit, and exposed nodes are exactly those $i=1,\ldots, n$ such that $x_i(t^*)=w_i$, $x_i(t^*)=0$, $0<x_i(^*)<w_i$, respectively, for the sequence $x(t)$ generated by the recursion \eqref{eq:iterative-algorithm} started in an arbitrary initial condition $x(0)$ in $\mc L_0^w$. Notice that, once such partition has been determined, in other to find all network equilibria, one is simply left to solve the linear system 
$$x_i=c_i+\sum_{j \ne i} P_{ji}x_j\,, \quad\; \forall\, i\in\mc V_0\,,$$
with boundary conditions $x_i=0$ for all $ i$ in $\mc V_-$ and 
$x_i=w_i$ for all $ i$ in $\mc V_+$, something that can be performed in finite time using standard algorithms for linear systems, e.g., Gaussian elimination. 
\end{remark}\medskip

\section{Geometry and uniqueness of network equilibria}\label{uniqueness}
In this section, we undertake a fundamental geometric study of the set of network equilibria and, in particular, we derive necessary and sufficient conditions for their uniqueness. We shall first consider two relevant special cases:
\begin{itemize}
	\item  when the matrix $P$ is asymptotically stable, i.e., such that $\rho(P)<1$ (Proposition \ref{prop:uniquenss-outconnected});
	\item when $P$ is irreducible and such that $\rho(P)=1$ (Proposition \ref{proposition:uniqueness-irreducible}).
\end{itemize}
Then, we build on these two cases in order to prove a general result (Theorem \ref{theorem:general}) on the geometric structure of the network equilibrium set $\mc X$ for every network $(P,w)$ such that $P$ has spectrum contained in the closed unitary disk. 


%
%

\begin{proposition}\label{prop:uniquenss-outconnected} 
For a network $(P,w)$ such that $\rho(P)<1$ and every exogenous flow $c$ in $\R^n$, there exists a unique network equilibrium $x$.
\end{proposition}
\proof
Let $x$ and $y$ in $\mc X$ be two network equilibria and put $\Delta=x- y$. We know from Corollary \ref{corollary:partition} (ii) that $\Delta_i=0$ for every $i$ in $\mc V_-\cup\mc V_+$. The proof is finished if $\mc V_0=\emptyset$. Otherwise, let 
 $z$ in $\R^{\mc V_0}$ and $Q$ in $\R^{\mc V_0\times\mc V_0}$ be the restrictions of $\Delta$ to $\mc V_0$ and of $P$ to $\mc V_0\times\mc V_0$, respectively. It then follows from Corollary \ref{corollary:partition} (i) that 
$z$ satisfies the equation $z=Q^{\top}z$. By Proposition \ref{proposition:PF} (ii), $\rho(Q)\le\rho(P)<1$, so that the matrix $(I-Q)$ is invertible and thus $z=0$. Therefore, $x=y$.
\qed

We now study the case of networks $(P,w)$ with $P$ irreducible and such that $\rho(P)=1$. 
The following result gives an explicit characterization of the condition of non-uniqueness as well as a representation of the set of network equilibria.

\begin{proposition}\label{proposition:uniqueness-irreducible}
Let $(P,w)$ be a network such that $P$ is irreducible and $\rho(P)=1$. Let  $\pi$ and $p$ be, respectively, left and right dominant eigenvectors of $P$, as in Proposition \ref{proposition:PF} (i).  
Then, for every exogenous flow $c$, there exists more than one network equilibrium in $\mc X$ if and only if
\be\label{eq:-min<min} p^{\top}c=0,\quad 	\min _{i} \l\{\frac{\nu_{i}}{\pi_{i}}\r\}+\min _{i}\left\{ \frac{w_{i}-\nu_{i}}{\pi_{i}}\right\}>0\,,\ee
where $\nu$ is any solution of the equation $\nu=P^{\top}\nu+c$ (c.f.~Proposition \ref{proposition:PF} (iv)). Moreover, in this case, the set of network equilibria is given by
\be\label{eq:segment}\mc X=\left\{x=\nu+\alpha\pi:\,-\min _{i} \l\{\frac{\nu_{i}}{\pi_{i}}\r\}\le\alpha\le\min _{i}\left\{ \frac{w_{i}-\nu_{i}}{\pi_{i}}\right\}\r\}\,.\ee
\end{proposition}
\begin{proof}
We first analyze the solution on $\R^n$ of the non-saturated linear system
\be\label{eq:linear-system}x=P^{\top}x+c\,.\ee  
Left multiplying by the vector $p$, we obtain
%
$$p^{\top}x=p^{\top}P^{\top}x+p^{\top}c=p^{\top}x+p^{\top}c$$ so that, for solutions of \eqref{eq:linear-system} to exist, it must hold  true that $p^{\top}c=0$.
On the other hand, if condition $p^{\top}c=0$ is satisfied, since $P$ is irreducible, Proposition \ref{proposition:PF} (iii) and (iv) ensure that 
the set of solutions of \eqref{eq:linear-system} is an affine line 
\be\label{line}\mc H=\{x=\nu+t\pi:\,t\in\R\}\,.\ee 
where $\nu$ is any solution of \eqref{eq:linear-system}.
Notice that solutions of the linear system \eqref{eq:linear-system} that belong to the complete lattice $\mc L_0^w$ are necessarily network equilibria, i.e., $\mc H \cap \mc L_0^w\subseteq \mc X\,.$
Moreover, observe that $\mc H\cap\mc L_0^w$ coincides with the right-hand side of \eqref{eq:segment} and that condition \eqref{eq:-min<min} is equivalent to saying that $\mc H\cap\mc L_0^w$ is a  segment of strictly positive length.

We are now ready to prove the statements of the theorem. 
Suppose first that there are multiple equilibria, i.e., $|\mc X|>1$. Since $P$ is irreducible, the only class of $\mc G_P$ is basic and final, so that Theorem \ref{theorem:partition} implies
that the node set partition \eqref{eq:partition-unique} is common to all network equilibria. If $\mc V_-\cup\mc V_+\neq\emptyset$, since $\mc V_0$ is a proper subset of $\mc V$, Proposition \ref{proposition:PF} (v) guarantees that the restriction $Q$ of $P$ to $\mc V_0\times\mc V_0$ has spectral radius smaller than $1$. Arguing exactly as in the proof of Proposition \ref{prop:uniquenss-outconnected}, we then deduce that $|\mc X|=1$ thus reaching a contradiction. Therefore, necessarily $\mc V=\mc V_0$. In this case, it follows from Corollary \ref{corollary:partition} (i) that all network equilibria are solutions of  \eqref{eq:linear-system}, i.e., $\mc H \cap \mc L_0^w= \mc X$. By our previous considerations, since this set is nonempty, the condition $p^{\top}c=0$ must hold true. Moreover, $|\mc X|>1$ implies that $\mc H \cap \mc L_0^w$ must be a segment of positive length that, as previously observed, is equivalent to the second condition in \eqref{eq:-min<min}. 

Suppose now that the conditions in \eqref{eq:-min<min} hold true. Then previous considerations imply that $\mc H \cap \mc L_0^w\subseteq\mc X$ is a segment of positive length. Non-uniqueness of network equilibria is thus proven.

Finally, notice that, if any of the two equivalent conditions hold, then $\mc H \cap \mc L_0^w= \mc X$ and this is equivalent to representation \eqref{eq:segment}. \qed
\end{proof}

\begin{remark}
The result above has a simple geometric interpretation in part already exploited in the proof. Assuming that $p^{\top}c=0$, the line $\mc H$ defined in \eqref{line} is the set of solutions of the non-saturated linear system \eqref{eq:linear-system}. The non-uniqueness condition \eqref{eq:-min<min} is simply the condition that this line intersects the interior part of the lattice $\mc L_0^w$ and the set of equilibria in this case is the segment obtained by this intersection. The minimal and maximal equilibria are the boundary points of this interval. We notice that the arguments used in the proof also show that, in the case of non-uniqueness, necessarily all nodes must be exposed nodes, namely $\mc V=\mc V_0$. 
\end{remark}

Below, we report an explicit calculation of the network equilibria for a three-dimensional network and two possible exogenous flows, respectively yielding uniqueness and multiplicity of network equilibria.

\medskip

\begin{example}\label{example:small-network} Consider the network $(P,w)$ where 
%
	$$P=\begin{bmatrix}
	0  &  0.75 &   0.25 \\
	0      &   0  &  1\\
	0.3 &    0.7 &         0
	\end{bmatrix}\,,\qquad w=\begin{bmatrix}5\\ 3\\ 2\end{bmatrix}\,.$$
	Notice that the matrix $P$ is stochastic and irreducible, hence we can take $p=\1$. The associated graph $\mc G_P$
	is depicted in Figure \ref{fig:example-small-network}. 	
\begin{figure}
		\centering	
		\begin{tikzpicture}[thick,scale=0.8, every node/.style={transform shape}]
		\draw  (-3.5,-1) ellipse (0.5 and 0.5);
		\draw  (0.5,-1) ellipse (0.5 and 0.5);
		\draw [-triangle 60] plot[smooth, tension=.7] coordinates {(0.5,-0.5) (0.5,0.5) (-1,1)};
		\draw [-triangle 60] plot[smooth, tension=.7] coordinates {(-3.5,-1.5) (-1.5,-2) (0.5,-1.5)};
		\node at (-3.5,-1) {1};
		\node at (0.5,-1) {2};
		\node at (-1.5,1) {3};
		\draw  (-1.5,1) ellipse (0.5 and 0.5);
		\draw  plot[smooth, tension=.7] coordinates {(-1,1)};
		\draw [-triangle 60] plot[smooth, tension=.7] coordinates {(-3.5,-0.5) (-3.5,0.5) (-2,1)};
		\draw [-triangle 60] plot[smooth, tension=.7] coordinates {(-1.5,0.5) (-2,-0.5) (-3,-1)};
		\draw [-triangle 60] plot[smooth, tension=.7] coordinates {(-1.5,0.5) (-1,-0.5) (0,-1)};
		\node at (-1.5,-2.5) {$ 0.75 $};
		\node at (-3.5,1) {$ 0.25 $};
		\node at (0.5,1) {$ 1 $};
		\node at (-2.5,-0.5) {$ 0.3 $};
		\node at (-0.5,-0.5) {$ 0.7 $};
		\end{tikzpicture}
		\caption{The network of Example  \ref{example:small-network}. \label{fig:example-small-network}}
		\label{f0N2}
	\end{figure}
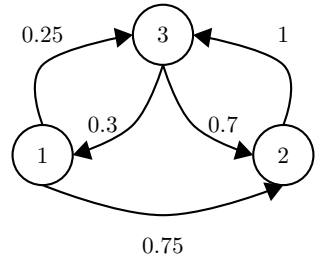
 We analyze uniqueness for two possible exogenous flows $$ c^{(1)}= [-1, 1, 0]^{\top}\,\qquad c^{(2)}= [-2, 2, 0]^{\top}\,.$$ 
 First of all, notice that
 $p^{\top}c^{(1)}=p^{\top}c^{(2)}=0$.
Moreover, a direct computation shows that
\begin{equation}\label{key}
\begin{array}{l}\ds\min _{i} \l\{\frac{\nu^1_{i}}{\pi_{i}}\r\}+\min _{i}\left\{ \frac{w_{i}-\nu^1_{i}}{\pi_{i}}\right\} \approx 1.60>0\\[10pt]
	\ds\min _{i} \l\{\frac{\nu^2_{i}}{\pi_{i}}\r\}+\min _{i}\left\{\frac{w_{i}-\nu^2_{i}}{\pi_{i}}\right\} \approx -6.41 <0.
	\end{array}
	\end{equation}
By Proposition \ref{proposition:uniqueness-irreducible} we deduce that for the flow $ c^{(1)}$ there are multiple equilibria, while for the flow $c^{(2)}$ the equilibrium is unique.
	The set of network equilibria $\mc X$ in the two cases is shown in Figure \ref{f3dd}. Notice how in the first case the line $\mc H$ has a non-trivial intersection with  the complete lattice $\mc L_0^w$ that is the segment of network equilibria. In contrast, in the second case, the line $\mc H$ does not intersect  the complete lattice $\mc L_0^w$ and the unique network equilibrium is a single point lying on the boundary of the lattice as some of its entries $x_i$ are necessarily saturated  at either $0$ or $w_i$.
\begin{figure}
		\centering
		\subfloat[The network $(P,w)$ with exogenous flow $c^{(1)}$ admits multiple equilibria (the black thick segment).]{\begin{tikzpicture}[thick,scale=0.8, x=0.75pt,y=0.75pt,yscale=-1,xscale=1]
	
	\draw  [color={rgb, 255:red, 69; green, 127; blue, 6 }  ,draw opacity=1 ][fill={rgb, 255:red, 126; green, 211; blue, 33 }  ,fill opacity=0.38 ] (277.53,221.02) -- (208.41,161.26) -- (209.5,84) -- (309.6,69.45) -- (378.72,129.21) -- (377.63,206.47) -- cycle ; \draw  [color={rgb, 255:red, 69; green, 127; blue, 6 }  ,draw opacity=1 ] (209.5,84) -- (278.62,143.76) -- (277.53,221.02) ; \draw  [color={rgb, 255:red, 69; green, 127; blue, 6 }  ,draw opacity=1 ] (278.62,143.76) -- (378.72,129.21) ;
	\draw    (159.3,211.2) -- (160.3,70.2) (155.37,201.17) -- (163.37,201.23)(155.44,191.17) -- (163.44,191.23)(155.51,181.17) -- (163.51,181.23)(155.58,171.17) -- (163.58,171.23)(155.65,161.17) -- (163.65,161.23)(155.73,151.17) -- (163.73,151.23)(155.8,141.17) -- (163.8,141.23)(155.87,131.17) -- (163.87,131.23)(155.94,121.17) -- (163.94,121.23)(156.01,111.17) -- (164.01,111.23)(156.08,101.17) -- (164.08,101.23)(156.15,91.17) -- (164.15,91.23)(156.22,81.17) -- (164.22,81.23)(156.29,71.18) -- (164.29,71.23) ;
	\draw    (159.3,211.2) -- (228.3,279.2) (169.23,215.37) -- (163.61,221.07)(176.35,222.39) -- (170.74,228.09)(183.48,229.41) -- (177.86,235.11)(190.6,236.43) -- (184.98,242.13)(197.72,243.45) -- (192.1,249.15)(204.84,250.47) -- (199.23,256.16)(211.97,257.49) -- (206.35,263.18)(219.09,264.51) -- (213.47,270.2)(226.21,271.52) -- (220.59,277.22) ;
	\draw    (228.3,279.2) -- (443.5,250) (237.67,273.89) -- (238.75,281.82)(247.58,272.55) -- (248.66,280.47)(257.49,271.2) -- (258.57,279.13)(267.4,269.86) -- (268.47,277.79)(277.31,268.51) -- (278.38,276.44)(287.22,267.17) -- (288.29,275.1)(297.13,265.82) -- (298.2,273.75)(307.04,264.48) -- (308.11,272.41)(316.94,263.14) -- (318.02,271.06)(326.85,261.79) -- (327.93,269.72)(336.76,260.45) -- (337.84,268.37)(346.67,259.1) -- (347.75,267.03)(356.58,257.76) -- (357.66,265.68)(366.49,256.41) -- (367.57,264.34)(376.4,255.07) -- (377.48,263)(386.31,253.72) -- (387.38,261.65)(396.22,252.38) -- (397.29,260.31)(406.13,251.03) -- (407.2,258.96)(416.04,249.69) -- (417.11,257.62)(425.95,248.35) -- (427.02,256.27)(435.86,247) -- (436.93,254.93) ;
	\draw [color={rgb, 255:red, 74; green, 144; blue, 226 }  ,draw opacity=1 ][line width=1.5]  [dash pattern={on 1.69pt off 2.76pt}]  (302,183.57) -- (395,227.57) ;
	\draw [color={rgb, 255:red, 74; green, 144; blue, 226 }  ,draw opacity=1 ][line width=1.5]  [dash pattern={on 1.69pt off 2.76pt}]  (196.5,126) -- (256.5,159) ;
	\draw  [color={rgb, 255:red, 0; green, 0; blue, 0 }  ,draw opacity=1 ][fill={rgb, 255:red, 0; green, 0; blue, 0 }  ,fill opacity=1 ] (298.75,183.57) .. controls (298.75,181.78) and (300.21,180.32) .. (302,180.32) .. controls (303.79,180.32) and (305.25,181.78) .. (305.25,183.57) .. controls (305.25,185.37) and (303.79,186.82) .. (302,186.82) .. controls (300.21,186.82) and (298.75,185.37) .. (298.75,183.57) -- cycle ;
	\draw [color={rgb, 255:red, 0; green, 0; blue, 0 }  ,draw opacity=1 ][line width=2.25]    (256.5,159) -- (302,183.57) ;
	\draw  [color={rgb, 255:red, 0; green, 0; blue, 0 }  ,draw opacity=1 ][fill={rgb, 255:red, 0; green, 0; blue, 0 }  ,fill opacity=1 ] (253.25,159) .. controls (253.25,157.21) and (254.71,155.75) .. (256.5,155.75) .. controls (258.29,155.75) and (259.75,157.21) .. (259.75,159) .. controls (259.75,160.79) and (258.29,162.25) .. (256.5,162.25) .. controls (254.71,162.25) and (253.25,160.79) .. (253.25,159) -- cycle ;
	
	\draw (120,146) node    {$x_{3}$};
	\draw (165,268) node    {$x_{1}$};
	\draw (335,301) node    {$x_{2}$};
	\draw (147,79) node    {$6$};
	\draw (146,104) node    {$4$};
	\draw (146,129) node    {$2$};
	\draw (146,153) node    {$0$};
	\draw (140,179) node    {$-2$};
	\draw (139,203) node    {$-4$};
	\draw (161,238) node    {$4$};
	\draw (195,267) node    {$2$};
	\draw (333,280) node    {$2$};
	\draw (391,272) node    {$0$};
	\draw (274,287) node    {$4$};
	\draw (218,290) node    {$0$};
	\draw (258,144) node    {$\overline{x}$};
	\draw (315,169) node    {$\underline{x}$};
	\draw (397,212) node  [color={rgb, 255:red, 74; green, 144; blue, 226 }  ,opacity=1 ]  {$\mathcal{H}$};
	\draw (299,107) node  [color={rgb, 255:red, 69; green, 125; blue, 3 }  ,opacity=1 ]  {$\mathcal{L}_{w}$};
	\draw (286,151.4) node [anchor=north west][inner sep=0.75pt]    {$\mathcal{X}$};

\end{tikzpicture}

}\quad
		\subfloat[The network $(P,w)$ with exogenous flow $c^{(2)}$ admits a unique equilibrium (the black dot).]{\begin{tikzpicture}[thick,scale=0.8, x=0.75pt,y=0.75pt,yscale=-1,xscale=1]

\draw  [color={rgb, 255:red, 69; green, 127; blue, 6 }  ,draw opacity=1 ][fill={rgb, 255:red, 126; green, 211; blue, 33 }  ,fill opacity=0.38 ] (277.53,221.02) -- (208.41,161.26) -- (209.5,84) -- (309.6,69.45) -- (378.72,129.21) -- (377.63,206.47) -- cycle ; \draw  [color={rgb, 255:red, 69; green, 127; blue, 6 }  ,draw opacity=1 ] (209.5,84) -- (278.62,143.76) -- (277.53,221.02) ; \draw  [color={rgb, 255:red, 69; green, 127; blue, 6 }  ,draw opacity=1 ] (278.62,143.76) -- (378.72,129.21) ;
\draw    (159.3,211.2) -- (160.3,70.2) (155.37,201.17) -- (163.37,201.23)(155.44,191.17) -- (163.44,191.23)(155.51,181.17) -- (163.51,181.23)(155.58,171.17) -- (163.58,171.23)(155.65,161.17) -- (163.65,161.23)(155.73,151.17) -- (163.73,151.23)(155.8,141.17) -- (163.8,141.23)(155.87,131.17) -- (163.87,131.23)(155.94,121.17) -- (163.94,121.23)(156.01,111.17) -- (164.01,111.23)(156.08,101.17) -- (164.08,101.23)(156.15,91.17) -- (164.15,91.23)(156.22,81.17) -- (164.22,81.23)(156.29,71.18) -- (164.29,71.23) ;
\draw    (159.3,211.2) -- (228.3,279.2) (169.23,215.37) -- (163.61,221.07)(176.35,222.39) -- (170.74,228.09)(183.48,229.41) -- (177.86,235.11)(190.6,236.43) -- (184.98,242.13)(197.72,243.45) -- (192.1,249.15)(204.84,250.47) -- (199.23,256.16)(211.97,257.49) -- (206.35,263.18)(219.09,264.51) -- (213.47,270.2)(226.21,271.52) -- (220.59,277.22) ;
\draw    (228.3,279.2) -- (443.5,250) (237.67,273.89) -- (238.75,281.82)(247.58,272.55) -- (248.66,280.47)(257.49,271.2) -- (258.57,279.13)(267.4,269.86) -- (268.47,277.79)(277.31,268.51) -- (278.38,276.44)(287.22,267.17) -- (288.29,275.1)(297.13,265.82) -- (298.2,273.75)(307.04,264.48) -- (308.11,272.41)(316.94,263.14) -- (318.02,271.06)(326.85,261.79) -- (327.93,269.72)(336.76,260.45) -- (337.84,268.37)(346.67,259.1) -- (347.75,267.03)(356.58,257.76) -- (357.66,265.68)(366.49,256.41) -- (367.57,264.34)(376.4,255.07) -- (377.48,263)(386.31,253.72) -- (387.38,261.65)(396.22,252.38) -- (397.29,260.31)(406.13,251.03) -- (407.2,258.96)(416.04,249.69) -- (417.11,257.62)(425.95,248.35) -- (427.02,256.27)(435.86,247) -- (436.93,254.93) ;
\draw [color={rgb, 255:red, 74; green, 144; blue, 226 }  ,draw opacity=1 ][line width=1.5]  [dash pattern={on 1.69pt off 2.76pt}]  (256.5,159) -- (395,227.57) ;
\draw [color={rgb, 255:red, 74; green, 144; blue, 226 }  ,draw opacity=1 ][line width=1.5]  [dash pattern={on 1.69pt off 2.76pt}]  (196.5,126) -- (256.5,159) ;
\draw  [color={rgb, 255:red, 0; green, 0; blue, 0 }  ,draw opacity=1 ][fill={rgb, 255:red, 0; green, 0; blue, 0 }  ,fill opacity=1 ] (275.25,161) .. controls (275.25,159.21) and (276.71,157.75) .. (278.5,157.75) .. controls (280.29,157.75) and (281.75,159.21) .. (281.75,161) .. controls (281.75,162.79) and (280.29,164.25) .. (278.5,164.25) .. controls (276.71,164.25) and (275.25,162.79) .. (275.25,161) -- cycle ;

\draw (120,146) node    {$x_{3}$};
\draw (165,268) node    {$x_{1}$};
\draw (335,301) node    {$x_{2}$};
\draw (147,79) node    {$6$};
\draw (146,104) node    {$4$};
\draw (146,129) node    {$2$};
\draw (146,153) node    {$0$};
\draw (140,179) node    {$-2$};
\draw (139,203) node    {$-4$};
\draw (161,238) node    {$4$};
\draw (195,267) node    {$2$};
\draw (333,280) node    {$2$};
\draw (391,272) node    {$0$};
\draw (274,287) node    {$4$};
\draw (218,290) node    {$0$};
\draw (397,212) node  [color={rgb, 255:red, 74; green, 144; blue, 226 }  ,opacity=1 ]  {$\mathcal{H}$};
\draw (299,107) node  [color={rgb, 255:red, 69; green, 125; blue, 3 }  ,opacity=1 ]  {$\mathcal{L}_{w}$};
\draw (288,147.4) node [anchor=north west][inner sep=0.75pt]    {$x$};

\end{tikzpicture}

		}
		\caption{\label{f3dd}  Sets of network equilibria for Example \ref{example:small-network}.}
	\end{figure}
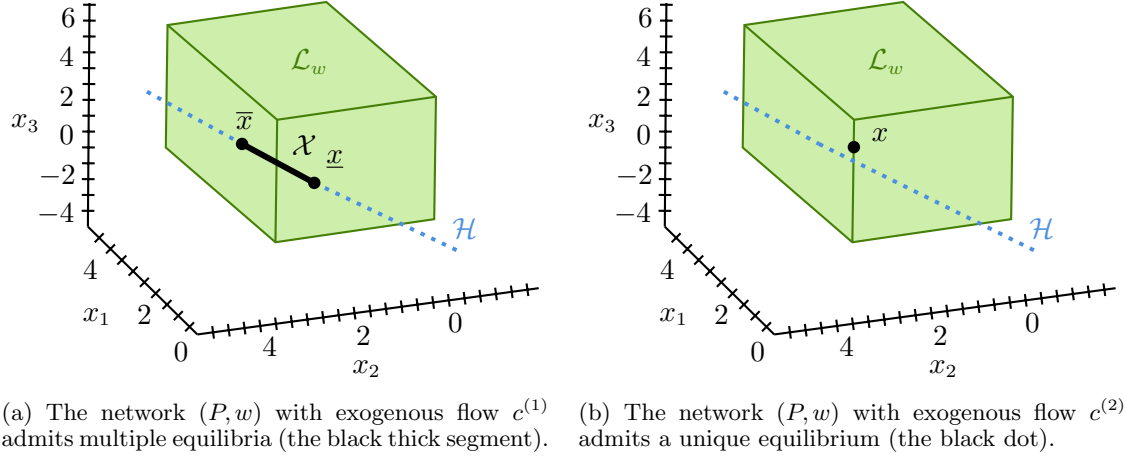	
	
\end{example}\medskip


  We now study the structure of network equilibria and give a full characterization of uniqueness in the general case of networks $(P,w)$ where $P$ is an arbitrary non-negative matrix with spectral radius $\rho(P)\leq 1$ and $w$ is an arbitrary non-negative vector. Our analysis relies on the partition of the node set in the classes of $P$ \begin{equation}\label{components}\mc V=\mc V_1\cup\cdots \cup \mc V_s\end{equation} 
and on the corresponding triangular structure of $P$ as described in \eqref{canonical}.

\begin{theorem}\label{theorem:general} Consider a network $(P,w)$ such that $\rho(P)\leq 1$, and an exogenous flow $c$.
Let \eqref{components} be the classes of $P$ and assume that $P$ is in the block triangular structure \eqref{canonical}. Indicate the related split of a vector $y$ in $\R^n$ as $y=[y^{(1)}, \dots , y^{(s)}]^{\top}$.
Then, the network equilibria $x$ in $\mc X$ iteratively satisfy the following properties:
	\begin{enumerate}
	\item[(i)] the projection $x^{(l)}$ on a class $\mc V_l$ such that $\rho(P^{(ll)})<1$ is unique;
	\item[(ii)] given $(x^{(1)},\dots , x^{(l-1)})$, the projection $x^{(l)}$ on a class $\mc V_l$ such that $\rho(P^{(ll)})=1$ is non-unique if and only if		
	\be\label{eq:necessary-nonunique}p^{(l)\top}\left(c^{(l)}+\sum_{1\le i< l}P^{(il)\top}x^{(i)}\right)=0\,,\ee
	and	\be\label{eq:min<min-3}\min_{i\in\mc V_l} \l\{\frac{\nu^{(l)}_i}{\pi^{(l)}_i}\r\}+\min_{i\in\mc V_l}\l\{\frac{w_i-\nu^{(l)}_i}{\pi^{(l)}_i}\r\}> 0\,,\ee
	where 
	\begin{itemize}
	\item $p^{(l)}=P^{(ll)}p^{(l)}$ is any right dominant eigenvector of the block $P^{(ll)}$;
	\item $\pi^{(l)}=P^{(ll)\top}\pi^{(l)}$ is any left dominant eigenvector of the block $P^{(ll)}$;
	\item $\nu^{(l)}=P^{(ll)\top}\nu^{(l)}+\sum_{i=1}^lP^{(jl)\top}x^{(i)}+ c^{(l)}$.
	\end{itemize}
	
	Moreover, in this case, given $[x^{(1)},\dots , x^{(l-1)}]^{\top}$, the projection $x^{(l)}$ of any equilibrium satisfies 
\be\label{eq:segment-2}x^{(l)}= \nu^{(l)}+\alpha\pi^{(l)}\,,\qquad -\min _{i\in\mc V_l} \l\{\frac{\nu^{(l)}_{i}}{\pi^{(l)}_{i}}\r\}\le\alpha\le\min _{i\in\mc V_l}\left\{ \frac{w_{i}-\nu^{(l)}_{i}}{\pi^{(l)}_{i}}\right\}\,.\ee

	\end{enumerate}
\end{theorem}
\proof It follows from \eqref{eq:nonlinear-system} and the block triangular structure of $P$ \eqref{canonical} that network equilibria satisfy the iterative relations
\begin{equation}\label{eq:nonlinear-system-iterative}
x^{(l)}=S_0^{w^{(l)}}\left(P^{(ll)\top}x^{(l)}+\sum_{0\le i< l}P^{(il)\top}x^{(i)}+c^{(l)}\right)\,,\qquad l=1,2,\ldots,s\,.
\end{equation}
The above says that the projection $x^{(l)}$ on the class $\mc V_l$ can be interpreted as a network equilibrium for the network $(P^{(ll)},w^{(l)})$ and exogenous flow $\sum_{i< l}P^{(il)\top}x^{(i)}+c^{(l)}$. The claim then follows from Propositions \ref{prop:uniquenss-outconnected} and \ref{proposition:uniqueness-irreducible}.
%
\qed

Notice that, as Proposition \ref{proposition:c} gives an efficient iterative way of computing the network equilibrium when this is unique, Theorem \ref{theorem:general} provides an explicit way of computing, in an iterative way, the entire lattice of network equilibria $\mc X$ in the general case when $\rho(P)\le1$.\medskip

\begin{remark}\label{remark:non-expansive} In the special case when the network is non-expansive (this includes the case when $P$ is stochastic or sub-stochastic) Theorem \ref{theorem:general} admits an important simplification. Indeed, in this case either $\rho(P)<1$, and then one can use Proposition \ref{prop:uniquenss-outconnected} directly to compute the unique network equilibrium (e.g., by using \eqref{eq:iterative-algorithm} as a distributed iterative algorithm, c.f. Remark \ref{remark:algorithm}), or $\rho(P)=1$ and the basic classes are final so that we can always assume that in the partition \eqref{components} they are the last ones. Precisely, in the latter case, we can assume that
\be\label{non-expansive-decomp}\rho(P^{(ll)})<1\; {\rm \: \: \:  for}\; \: \: \:  l\leq m,\quad  \quad  \rho(P^{(ll)})=1\; {\rm \: \: \: for}\; \: \: \:  m<l\leq  s\,.\ee
The projection $(x^{(1)},\dots , x^{(m)})$ of the network equilibria $x$ on the first $m$ classes is unique. For each basic class $\mc V_l$, with $m<l\le s$, the non uniqueness condition of the projection $x^{(l)}$ is given by
\be\label{non-expansive-decomp2}p^{(l)\top}\left(c^{(l)}+\sum_{1\le i\le m}P^{(il)\top}x^{(i)}\right)=0\,,\ee
together with \eqref{eq:min<min-3}. We notice that these conditions only depend on $(x^{(1)},\dots , x^{(m)})$. In other words, once the solution on the non-basic classes is computed, the check of uniqueness and the parametrization of the solutions in case of non-uniqueness in the various basic classes are completely decoupled. 
\end{remark}\medskip

 \begin{remark}  \label{remark:bifurcation}
Notice that our analysis has mostly focused on networks $(P,w)$ with spectral radius $\rho(P)\le1$. In fact, Theorem \ref{theorem:general} provides a complete description of the set of network equilibria $\mc X$ in this case. It is worth stressing out that, for networks with $\rho(P)>1$, while $\mc X$ remains a nonempty complete lattice as per Proposition \ref{proposition:c}, its geometry can differ quite significantly in this case. In fact, consider a simple example with a single node, $P=2$, and $w=1$. Then, depending of the value of the exogenous flow $c$ in $\R$ the set of network equilibria is 
\be\label{eq:bifurcation}\mc X=\left\{\ba{lcl}\{0\}&\se& c<-1\\[7pt]\{0,-c,1\}&\se&-1\le c\le0\\[7pt]\{1\}&\se&c>0\,,\ea\right.\ee
as illustrated in Figure \ref{fig:expansive-X}. 
In particular, notice that for values of the exogenous flow $c$ in $\mc M=[-1,0]$, there are multiple isolated network equilibria, specifically $|\mc X|=2$ for $c=-1$ and  $|\mc X|=3$ for $-1<c<0$. This is in stark contrast with the case $\rho(P)\le1$, where Theorem \ref{theorem:general} in particular implies that, when the network equilibrium is not unique, there is in fact a continuum of network equilibria. 
\end{remark}

\begin{figure}\begin{center}
\includegraphics[height=4cm]{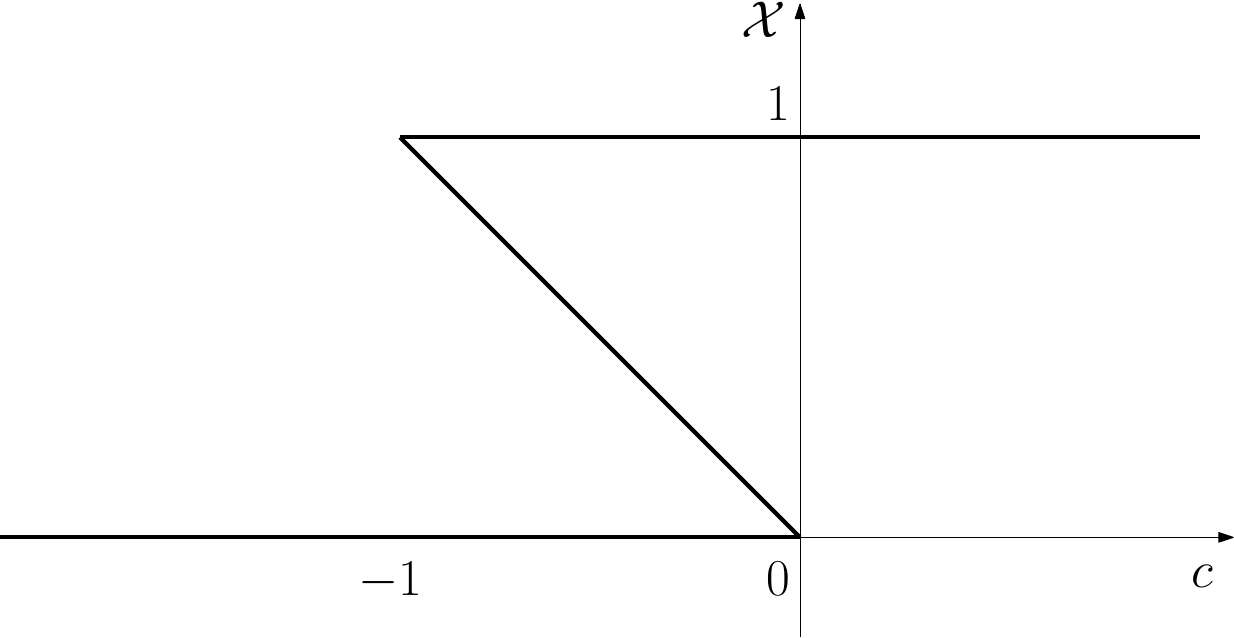}
\end{center}
\caption{The set of network equilibria $\mc X$ for the network discussed in Remark \ref{remark:bifurcation} as a function of the exogenous flow $c$.\label{fig:expansive-X}}
\end{figure}

\section{Continuity of network equilibria and the lack thereof}\label{discontinuous}
 
In this section, we study the dependence of the network equilibria of a given network $(P,w)$ on the exogenous flow $c$. This analysis is crucial to study the way exogenous shocks affect the payment equilibria in financial networks (c.f.~Section \ref{sec:financial}) or the individual marginal benefits affect the Nash equilibrium in quadratic network games (c.f.~ Section \ref{sec:game}). 


Let us consider a given network $(P,w)$ and use the notation
$$\mc X(c),\quad \ov x(c),\quad \underline x(c)$$
to emphasize the dependence of, respectively, the set of network equilibria, and the maximal and minimal network equilibrium on the exogenous flow $c$.
%
%
%
%
Moreover, let
\be\label{eq:UM-def}\mc U=\left\{c\in \R^n:\,|\mc X(c)|=1\right\}\,,\qquad \mc M=\R^n\setminus\mc  U\,,\ee
be the subsets of exogenous flows for which the network equilibrium is unique and, respectively, there are multiple  network equilibria.  For exogenous flows $c$ in $\mc U$, we shall also use the notation $$x(c)=\ul x(c)=\ov x(c)$$  for the unique equilibrium.

The following result gives a complete picture of the behavior of the set of network equilibria $\mc X(c)$ as a function of the exogenous flow $c$. It shows that
the set of exogenous flows $\mc M$ for which the network equilibrium is not unique  has Lebesgue measure $0$ and is contained in the union of a finite number of graphs of continuous functions. Moreover,  the network equilibrium $x(c)$ is a piece-wise continuous function of the exogenous flow $c$ that undergoes jump discontinuities when $c$ crosses the non-uniqueness set $\mc M$. 
\begin{theorem}\label{theo:continuity} For a network $(P,w)$ such that $\rho(P)\le 1$, let $m$ be number of basic classes of $P$ and 
let $\mc U$ and $\mc M$ be defined as in \eqref{eq:UM-def}. 
Then, 
\begin{enumerate}
\item[(i)]  the non-uniqueness set $\mc  M$ has  Lebesgue measure $0$ and is contained in the closed set consisting of the union of at most $m$ graphs of scalar continuous functions;  
%
%
\item[(ii)] the map $c\mapsto x(c)$ is continuous on the uniqueness set $\mc U$; 
\item[(iii)] for every exogenous flow $c^*$ in $\mc M$, 
$$ \liminf\limits_{{\substack{c\in\mc U\\ c\to c^*}}} x(c)=\underline x(c^*)\,,\qquad \limsup\limits_{\substack{c\in\mc U\\ c\to c^*}}x(c)=\ov x(c^*)\,.$$
\end{enumerate}
\end{theorem}

\proof We start with a preliminary computation that will  prove useful in the following derivations. Consider a sequence $c(1),c(2),\ldots$ of exogenous flows in $\R^n$ such that 
\be\label{sequence}c{(t)}\stackrel{t\to+\infty}{\longrightarrow} c^*\,,\qquad\underline x(c{(t)})\stackrel{k\to+\infty}{\longrightarrow} x^*\,.\ee Since $$\underline x(c(t))=S_0^w\left(P^{\top}\ul x(c(t))+c(t)\right)\,,$$ 
for all $t=1,2,\ldots$, passing to the limit in both sides of the above, by continuity we get that 
$$x^*=S_0^w\left(P^{\top} x^*+c^*\right)\,,$$ thus showing that  $x^*$ belongs to $\mc X(c^*)$. In particular, this implies that  
\be\label{limit}\underline x(c^*)\leq x^*\leq \ov x(c^*)\,.\ee
Arbitrariness of the sequence satisfying \eqref{sequence} and \eqref{limit} imply that 
\be\label{liminfsup}  \underline x(c^*)\leq \liminf\limits_{{\substack{c\in\mc U\\ c\to c^*}}} x(c)\leq\limsup\limits_{{\substack{c\in\mc U\\ c\to c^*}}} x(c)\leq \ov x(c^*)\ee
In particular, for every exogenous flow $c^*$ in $\mc U$, we have that $\underline x(c^*)=\ov x(c^*)$ and then relation \eqref{liminfsup} yields point (ii) of the claim. 

Consider now the partition \eqref{components} of the node set into the classes of $P$ and assume without loss of generality that $P$ is in the block triangular structure \eqref{canonical}. As usual, we indicate the relative split of any vector $y$ in $\R^n$ as $y=[y^{(1)}, \dots , y^{(s)}]^{\top}$. Assume that $l_1<\cdots <l_m$ are the indices among $\{1,\dots , s\}$ corresponding to the basic classes $\mc V_{l_1},\dots , \mc V_{l_m}$. For a fixed $j$, we consider the projection of the set of equilibria on $\mc V_1\cup\cdots\cup \mc V_{l_j-1}$. Notice that, because of the triangular structure of $P$, such projected set depends on $c=[c^{(1)},\dots c^{(s)}]^{\top}$ only through the sub-vector $
[c^{(1)},\dots c^{(l_j-1)}]^{\top}$. Suppose that for a given $c$ and for a given $j$, such projected set is a singleton and indicate the projected block components of such equilibrium as $x^{(i)}([c^{(1)},\dots c^{(l_j-1)}])$ for $i=1,\dots ,l_j-1$.
It then follows from Theorem \ref{theorem:general} that a necessary condition for the projection of the equilibria on $\mc V_{(l_j)}$ not to be unique, is that 
\begin{equation}\label{unique}p^{(l_j)\top}\left(c^{(l_j)}+\sum_{i< l_j}P^{(il_j)\top}x^{(i)}[c^{(1)},\dots c^{(l_j-1)}]\right)=0\end{equation}
Now, define the sets $\mc U_k, \mc M_k\subseteq \R^{\mc V_1\cup\cdots\cup \mc V_{l_k}}$ as follows: 
\begin{equation}\label{unique2}\begin{array}{rcl} 
\mc U_k&=&\{[c^{(1)}, \dots , c^{(l_k)}]:\, [c^{(1)},\dots c^{(l_j)}]\, \mbox{\rm does not satisfy \eqref{unique} }\; \forall j\leq k
\}\,;\\[7pt]
\mc M_k&=&\{[c^{(1)}, \dots , c^{(l_k)}]:\, [c^{(1)},\dots c^{(l_j)}]\in \mc U_j\; \forall j\leq k-1 , \mbox{\rm and \eqref{unique} is satisfied for}\;  j= k
\}\,.
\end{array}
\end{equation}
%
%
%
%
%
%
%
%
%
Put $\tilde{\mc M}_k=\mc M_k\times \R^{\mc V_{l_k+1}\cup\cdots\cup \mc V_{m}}$ and notice that the
considerations above imply that \begin{equation}\label{closed-nonunique}\mc M\subseteq \bigcup_{k=1}^m\tilde{\mc M}_k\end{equation}
Applying item (ii) to the restricted network consisting of the nodes in $\mc V_1\cup\cdots\cup \mc V_{l_k}$ we deduce that, for every $i=1,\dots ,l_k$, the functions $x^{(i)}([c^{(1)},\dots c^{(l_k)}])$ are continuous on the set $\mc U_k$. This fact, together with the definition of $\mc M_k$ and the form of condition \eqref{unique}, allows us to conclude that $\mc M_k$ is the graph of a continuous function defined on $\mc U_{k-1}\times \R^{\mc V_{l_k}\setminus\{s_k\} }$ where $s_k$ is any element in $\mc V_{l_k}$. An analogous conclusion then holds true for $\tilde{\mc M}_k$. This proves (i).


We are now left with proving (iii). Let $c^*$ in $\mc M$ be an exogenous flow giving rise to multiple equilibria and define the sequence of exogenous flows $c(t)$ as follows:
$$c^{(i)}(t)=c^{*(i)}-\frac{1}{t}p^{(i)}\; \: \: \:  \forall \:  i=1,\dots , s\,.$$
where $p^{(i)}$ is any right dominant eigenvector of the block $P^{(ii)}$
We claim that $c(t)$ necessarily belongs to $\mc U$ for  sufficiently large $t$. Indeed, a simple iterative argument shows that, if $t$ is sufficiently large,  $[c^{(1)}(t),\dots c^{(l_k)}(t)]\in\mc U_k$ for every $k$ and therefore $c(t)\not\in\tilde{\mc M}_k$ for every $k$. The claim then follows from \eqref{closed-nonunique}. 
%
%
Since $c(t)\le c^*$ for every $t=1,2,\ldots$, it follows from Proposition \ref{proposition:c} (iv) that
$$ x(c(t))=\underline x(c(t))\leq \underline x(c^*)\,.$$
Using relation \eqref{liminfsup}, we deduce that
\be\label{liminfeq}  \liminf\limits_{{\substack{c\in\mc U\\ c\to c^*}}}  x(c)=\underline x(c^*)\,.\ee
An analogous argument allows us to prove the other relation in (iii) concerning the $\limsup$.
%
%
%
%
%
%
%
%
%
%
%
%
%
%
\qed\medskip


For the special case of non-expansive networks $(P,w)$, we are able to characterize the maximum discontinuity jump of the network equilibrium as the exogenous flow $c$ varies in $\R^n$, as stated in the following result.

\begin{corollary}\label{coro:jump} For a non-expansive network $(P,w)$, consider the partition \eqref{components} of the node set into the classes of $P$ and let the block triangular structure of $P$ be as in \eqref{canonical}. Let $\pi^{(l)}$ be any left dominant eigenvalue relative to $P^{(ll)}$. Then,
\begin{enumerate}
\item[(i)] for every exogenous flow $c$, indicated with 
$$L_c=\{l=1,\dots , s\,|\, \mc V_l\;\hbox{\rm is basic and}\; (\ref{non-expansive-decomp2})\, \hbox{\rm is satisfied}\}$$ the norm of the jump discontinuity of the network equilibrium at $c$ can be expressed as 
\be\label{eq:jump}\|\ov x(c)-\ul x(c)\|_p^p=\sum_{\substack{l=1,\ldots,s:\\ l\in L_c}}\left(\left[\min_{i\in\mc V_l}\frac{w_i-\nu_i^{(l)}}{\pi_i^{(l)}}+\min_{i\in\mc V_l}\frac{\nu^{(l)}_i}{\pi_i^{(l)}}
\right]^+\right)^p\|\pi^{(l)}\|_p^p\,,\ee 
where $\nu^{(l)}$ is defined in Theorem \ref{theorem:general}.
\item[(ii)]
the maximum jump discontinuity norm is for $c=0$ and is given by
\be\label{eq:max-jump}\max_{c\in\R^n}\|\ov x(c)-\ul x(c)\|_p^p=\|\ov x(0)-\ul x(0)\|_p^p=\sum_{\substack{l=1,\ldots,s:\\ \mc V_l \,{\rm basic}}}\left(\min_{i\in\mc V_l}\frac{w_i}{\pi_i^{(l)}}\right)^p\|\pi^{(l)}\|_p^p\,,\ee 
\end{enumerate}
%
%
\end{corollary}
\proof 
Formula \eqref{eq:jump} directly follows from Theorem \ref{theorem:general} by virtue of the non-uniqueness condition \eqref{eq:necessary-nonunique} as modified  in \eqref{non-expansive-decomp2} and the structure of solutions as expressed in \eqref{eq:segment-2}.
From \eqref{eq:jump}, we obtain that
$$\|\ov x(c)-\ul x(c)\|_p^p\leq \sum_{\substack{l=1,\ldots,s:\\ l\in L_c}}\left(\min_{i\in\mc V_l}\frac{w_i}{\pi_i^{(l)}}\right)^p\|\pi^{(l)}\|_p^p\leq \sum_{\substack{l=1,\ldots,s:\\ \mc V_l \,{\rm basic}}}\left(\min_{i\in\mc V_l}\frac{w_i}{\pi_i^{(l)}}\right)^p\|\pi^{(l)}\|_p^p$$
On the other hand, since for $c=0$ every $l$ for which $\mc V_l$ is a basic class belongs to $L_c$, and since we can choose $\nu^{(l)}=0$,  formula \eqref{eq:jump} yields \eqref{eq:max-jump}.

 \qed

A few comments are in order. 
First, notice that, for networks such that $\rho(P)=1$, Theorems \ref{theorem:general} and  \ref{theo:continuity} ensure that the network equilibrium is generically unique and at the same time characterize the set $\mc M$ of exogenous flows inducing multiple network equilibria.  As a function of the exogenous flow $c$, the network equilibrium $x(c)$ is proven to be a piece-wise continuous function (it is also monotone in $c$ thanks to Proposition \ref{proposition:c}) with jump discontinuities occurring exactly when crossing the non-uniqueness set $\mc M$. For the relevant family of non-expansive networks, Corollary \ref{coro:jump} establishes an explicit formula for the value norm of these jumps. For networks with $\rho(P)<1$, Proposition \ref{proposition:uniqueness-irreducible} guarantees that the network equilibrium $x(c)$ is unique for every value of the exogenous flow $c$ and, in this case, it is a monotone continuous function of it. 

%

{ 
Another relevant observation is that the multiplicity of network equilibria for networks $(P,w)$ with spectral radius $\rho(P)=1$ and particular exogenous flows $c^*$ can also be interpreted as an indicator of high sensitivity in the dependence of the network equilibrium $\tilde x(c)$ of networks $(\tilde P,w)$ with spectral radius $\rho(\tilde P)<1$ that are sufficiently close to the nominal network $(P,w)$. This is first illustrated by the following simple example. \medskip

\begin{example}\label{example:blowup} 
Consider the family of networks $(P^{(\eps)},w)$, indexed by $\eps \in [0,1)$, with $n=2$ nodes and 
$$P^{(\eps)}=\left[\ba{cc}1-\eps&1\\0&1/2\ea\right]\,,\qquad w=\left[\ba{c}2\\1\ea\right]\,.$$
Notice that for $\eps \in (0,1)$ we have $\rho(P^{(\eps)})=\max\{1-\eps,1/2\}$ and for every exogenous flow $c$ in $\R^2$ there exists a unique network equilibrium $x^{(\eps)}(c)$ with entries 
$$x_1^{(\eps)}(c)=S_0^2(c_1/\eps)\,,\qquad x_2^{(\eps)}(c)=S_0^1(2c_2+2S_0^2(c_1/\eps))\,.$$
On the other hand, for $\eps=0$ we recover the same network as in Example \ref{example:non-partition}, with $\rho(P^{(0)})=1$. For such network, the set of exogenous flows giving rise to multiple equilibria is the whole line $\mc M=\{(0,t):\,t\in\R\}$. It is then clear as the sensitivity of the first entry of the network equilibrium satisfies 
$$\frac{\partial x_1^{\eps}}{\partial c_1}(0^+,c_2)=\frac1{\eps}\stackrel{\eps\downarrow0}{\longrightarrow}+\infty\,,$$
for every $c_1$ in $\R$. 
\end{example}\medskip

We conclude this section by discussing implications of our results in the two main motivating applications presented in Section \ref{sec:motivations}. 

\subsection{Systemic risk in financial networks}\label{sec:systemicrisk}
Consider the generalized Eisenberg and Noe financial network model introduced in Section \ref{sec:financial}. 
In order to measure the aggregated effect of a shock, it is useful to introduce a risk measure known as \emph{systemic loss}  \cite{Peyton2002}.  Let $c^{\circ}$ be a nominal exogenous flow for which all nodes in the financial network are fully liable, i.e., such that $x(c^{\circ})=w$. The, let $c\le c^{\circ}$  be the exogenous flow after a shock has negatively affected the assets and external credits of some of the financial entities in the network and let $x(c)$ be a corresponding network equilibrium. As in Subsection \ref{sec:financial}, let the net worth vectors before and after the shock be, respectively,
$v^{\circ}=P^{\top}w+c^{\circ}-w$ and $ v=P^{\top}x(c)+c-w$. Then, the systemic loss is defined as their aggregate difference 
%
%
%
  \begin{equation}\label{systemic-loss}
  	l\left(c^{\circ}, c\right) :=\1^{\top}\left(v^{\circ}-v\right)=\1^{^{\top}}\left(P^{^{\top}} w+c^{\circ}-w-\left(P^{^{\top}} x(c)+c-w\right)\right) =\1^{^{\top}}(c^{\circ}-c)+\1^{\top}(w-x(c))\,.
  \end{equation}
  In the rightmost side of the expression above, the term $\1^{\top}(c^{\circ}-c)$ represents the direct loss inflicted by the shock, while $\1^{\top}(w-x(c))$ represents the indirect loss triggered by reduced payments and is also referred to as shortfall term.  Then, we may apply \eqref{systemic-loss} and Theorem \ref{theo:continuity} (iii) to obtain the following expression for the size of the jump discontinuity of the systemic loss at some point $c=c^*$:
 \begin{equation}\label{eq:exdisc}
  		\Delta l\left(c^{*}\right) :=\limsup _{c \in \mathcal{U} \atop c \rightarrow c^{*}} l\left(c^{\circ}, c\right)-\liminf _{c \in \mathcal{U} \atop c \rightarrow c^{*}} l\left(c^{\circ}, c\right)
  		=\|\ov x(c^*)-\ul x(c^*)\|_1\,.
%
  \end{equation}
Explicit estimates of the expression above can then be obtained using formula (\ref{eq:jump}) in Corollary \ref{coro:jump}. Systemic loss jumps are expected to play a crucial role in the resilience analysis of the financial network as they will often be associated
to important failure events where several nodes simultaneously lose  their liability, as illustrated in the example below. 

\medskip


 \begin{example}\label{example:loss}
%
%
%
%
Consider the financial network of Example \ref{example:small-network}. 
(Figure \ref{f0N2}). 
The set $\mc M$ of exogenous flows giving rise to multiple network equilibria is plotted in Figure \ref{fig: setm}. 
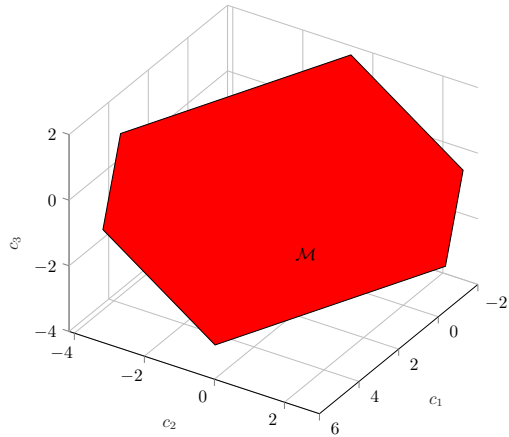
\begin{figure}
	\centering
\begin{tikzpicture}[scale=0.6]
	
	\begin{axis}[%
		width=3.557in,
		height=3.566in,
		at={(0.597in,0.481in)},
		scale only axis,
		xmin=-2,
		xmax=6,
		tick align=outside,
		xlabel style={font=\color{white!15!black}},
		xlabel={$c_1$},
		ymin=-4.15,
		ymax=3,
		ylabel style={font=\color{white!15!black}},
		ylabel={$c_2$},
		zmin=-4,
		zmax=2,
		zlabel style={font=\color{white!15!black}},
		zlabel={$c_3$},
		view={122.399999034916}{37.7999999577879},
		axis background/.style={fill=white},
		axis x line*=bottom,
		axis y line*=left,
		axis z line*=left,
		xmajorgrids,
		ymajorgrids,
		zmajorgrids,
		legend style={at={(1.03,1)}, anchor=north west, legend cell align=left, align=left, draw=white!15!black}
		]
		
		\addplot3[area legend, draw=black, fill=red]
		table[row sep=crcr] {%
			x	y	z\\
			-0.36	3	-2.64\\
			-1.96	2.6	-0.640000000000001\\
			-1.6	-0.4	2\\
			3.4	-4.15	0.749999999999998\\
			5	-3.75	-1.25\\
			4.64	-0.75	-3.89\\
		}--cycle;
		
	\end{axis}
	
	\begin{axis}[%
		width=5.833in,
		height=4.375in,
		at={(0in,0in)},
		scale only axis,
		xmin=0,
		xmax=1,
		ymin=0,
		ymax=1,
		axis line style={draw=none},
		ticks=none,
		axis x line*=bottom,
		axis y line*=left,
		legend style={legend cell align=left, align=left, draw=white!15!black}
		]
		\node[below right, align=left, draw=red]
		at (rel axis cs:0.432,0.452) {$\mc M$};
	\end{axis}
\end{tikzpicture}%
\caption{The set of critical shocks $\mc M$.}
\label{fig: setm}
\end{figure}
   Consider an initial exogenous flow $c^{\circ}=[5, 2, 2]^{\top}$ and a perturbation of it 
$c=c^{\circ}-\epsilon q$, where $q= [
0.07,
0.59,
0.34]^{\top}$, and $\epsilon \in [0,14]$. 
A straightforward computation, using condition of Proposition \ref{proposition:uniqueness-irreducible}, implies that the only case where we have multiple equilibria is for $ \epsilon=9$ corresponding to the exogenous flow $c^*=[4.4, -3.3, -1.1]^{\top}$
 for which$$\Delta l(c^*)=\min _{i}\left\{\nu_{i}/\pi_{i}\right\}+\min _{i}\left\{(w_{i}-\nu_{i})/{\pi_{i}}\right\} \approx 4.44>0$$

The loss function and the equilibrium $x$ as functions of $\epsilon$ are plotted in Figure \ref{lossdisc}.
\begin{figure}
	\centering
	\subfloat[Loss as a function of $\epsilon$. We can see the jump discontinuity at $\epsilon=9$. ]
	{\begin{tikzpicture}[thick,scale=0.97, x=0.75pt,y=0.75pt,yscale=-1,xscale=1]

\draw [line width=0.75]  (71,241.49) -- (325,241.49)(96.4,74.89) -- (96.4,260) (318,236.49) -- (325,241.49) -- (318,246.49) (91.4,81.89) -- (96.4,74.89) -- (101.4,81.89) (126.4,236.49) -- (126.4,246.49)(156.4,236.49) -- (156.4,246.49)(186.4,236.49) -- (186.4,246.49)(216.4,236.49) -- (216.4,246.49)(246.4,236.49) -- (246.4,246.49)(276.4,236.49) -- (276.4,246.49)(306.4,236.49) -- (306.4,246.49)(91.4,211.49) -- (101.4,211.49)(91.4,181.49) -- (101.4,181.49)(91.4,151.49) -- (101.4,151.49)(91.4,121.49) -- (101.4,121.49)(91.4,91.49) -- (101.4,91.49) ;
\draw   (133.4,253.49) node[anchor=east, scale=0.75]{2} (163.4,253.49) node[anchor=east, scale=0.75]{4} (193.4,253.49) node[anchor=east, scale=0.75]{6} (223.4,253.49) node[anchor=east, scale=0.75]{8} (253.4,253.49) node[anchor=east, scale=0.75]{10} (283.4,253.49) node[anchor=east, scale=0.75]{12} (313.4,253.49) node[anchor=east, scale=0.75]{14} (93.4,211.49) node[anchor=east, scale=0.75]{5} (93.4,181.49) node[anchor=east, scale=0.75]{10} (93.4,151.49) node[anchor=east, scale=0.75]{15} (93.4,121.49) node[anchor=east, scale=0.75]{20} (93.4,91.49) node[anchor=east, scale=0.75]{25} ;
\draw  [draw opacity=0] (95.4,91.89) -- (314,91.89) -- (314,240.85) -- (95.4,240.85) -- cycle ; \draw  [color={rgb, 255:red, 0; green, 0; blue, 0 }  ,draw opacity=0.11 ] (95.4,91.89) -- (95.4,240.85)(125.4,91.89) -- (125.4,240.85)(155.4,91.89) -- (155.4,240.85)(185.4,91.89) -- (185.4,240.85)(215.4,91.89) -- (215.4,240.85)(245.4,91.89) -- (245.4,240.85)(275.4,91.89) -- (275.4,240.85)(305.4,91.89) -- (305.4,240.85) ; \draw  [color={rgb, 255:red, 0; green, 0; blue, 0 }  ,draw opacity=0.11 ] (95.4,91.89) -- (314,91.89)(95.4,121.89) -- (314,121.89)(95.4,151.89) -- (314,151.89)(95.4,181.89) -- (314,181.89)(95.4,211.89) -- (314,211.89) ; \draw  [color={rgb, 255:red, 0; green, 0; blue, 0 }  ,draw opacity=0.11 ]  ;
\draw [color={rgb, 255:red, 74; green, 144; blue, 226 }  ,draw opacity=1 ][line width=1.5]    (96.4,241.49) -- (199,204) ;

\draw [color={rgb, 255:red, 74; green, 144; blue, 226 }  ,draw opacity=1 ][line width=1.5]    (199,204) -- (230,182) ;

\draw [color={rgb, 255:red, 255; green, 0; blue, 31 }  ,draw opacity=1 ][line width=1.5]    (230,154) -- (230,182) ;

\draw [color={rgb, 255:red, 74; green, 144; blue, 226 }  ,draw opacity=1 ][line width=1.5]    (230,154) -- (303,100) ;

\draw [color={rgb, 255:red, 255; green, 0; blue, 0 }  ,draw opacity=1 ] [dash pattern={on 4.5pt off 4.5pt}]  (230,182) -- (230,240) ;

\draw  [color={rgb, 255:red, 254; green, 0; blue, 0 }  ,draw opacity=1 ][fill={rgb, 255:red, 255; green, 0; blue, 0 }  ,fill opacity=1 ] (228.44,241.56) .. controls (228.44,240.7) and (229.14,240) .. (230,240) .. controls (230.86,240) and (231.56,240.7) .. (231.56,241.56) .. controls (231.56,242.42) and (230.86,243.11) .. (230,243.11) .. controls (229.14,243.11) and (228.44,242.42) .. (228.44,241.56) -- cycle ;

\draw (337,238.89) node   {$\epsilon $};
\draw (66,77.89) node   {$l( \epsilon )$};
\draw (91,248.49) node [scale=0.7]  {$0$};
\draw (230,254) node [scale=0.7,color={rgb, 255:red, 255; green, 0; blue, 0 }  ,opacity=1 ]  {$9$};

\end{tikzpicture}

	} \quad
	\subfloat[Solution vector as a function of $\epsilon$. We can see the jump discontinuity at $\epsilon=9$, which brings nodes 1 and 3 to suddenly default.]
	{\begin{tikzpicture}[thick,scale=0.97, x=0.75pt,y=0.75pt,yscale=-1,xscale=1]
	
	\draw [line width=0.75]  (71,241.49) -- (325,241.49)(96.4,74.89) -- (96.4,260) (318,236.49) -- (325,241.49) -- (318,246.49) (91.4,81.89) -- (96.4,74.89) -- (101.4,81.89) (126.4,236.49) -- (126.4,246.49)(156.4,236.49) -- (156.4,246.49)(186.4,236.49) -- (186.4,246.49)(216.4,236.49) -- (216.4,246.49)(246.4,236.49) -- (246.4,246.49)(276.4,236.49) -- (276.4,246.49)(306.4,236.49) -- (306.4,246.49)(91.4,211.49) -- (101.4,211.49)(91.4,181.49) -- (101.4,181.49)(91.4,151.49) -- (101.4,151.49)(91.4,121.49) -- (101.4,121.49)(91.4,91.49) -- (101.4,91.49) ;
\draw   (133.4,253.49) node[anchor=east, scale=0.75]{2} (163.4,253.49) node[anchor=east, scale=0.75]{4} (193.4,253.49) node[anchor=east, scale=0.75]{6} (223.4,253.49) node[anchor=east, scale=0.75]{8} (253.4,253.49) node[anchor=east, scale=0.75]{10} (283.4,253.49) node[anchor=east, scale=0.75]{12} (313.4,253.49) node[anchor=east, scale=0.75]{14} (93.4,211.49) node[anchor=east, scale=0.75]{1} (93.4,181.49) node[anchor=east, scale=0.75]{$w_3=$2} (93.4,151.49) node[anchor=east, scale=0.75]{$w_2=$3} (93.4,121.49) node[anchor=east, scale=0.75]{4} (93.4,91.49) node[anchor=east, scale=0.75]{$w_1=$5} ;
	\draw  [draw opacity=0] (95.4,91.89) -- (314,91.89) -- (314,240.85) -- (95.4,240.85) -- cycle ; \draw  [color={rgb, 255:red, 0; green, 0; blue, 0 }  ,draw opacity=0.11 ] (95.4,91.89) -- (95.4,240.85)(125.4,91.89) -- (125.4,240.85)(155.4,91.89) -- (155.4,240.85)(185.4,91.89) -- (185.4,240.85)(215.4,91.89) -- (215.4,240.85)(245.4,91.89) -- (245.4,240.85)(275.4,91.89) -- (275.4,240.85)(305.4,91.89) -- (305.4,240.85) ; \draw  [color={rgb, 255:red, 0; green, 0; blue, 0 }  ,draw opacity=0.11 ] (95.4,91.89) -- (314,91.89)(95.4,121.89) -- (314,121.89)(95.4,151.89) -- (314,151.89)(95.4,181.89) -- (314,181.89)(95.4,211.89) -- (314,211.89) ; \draw  [color={rgb, 255:red, 0; green, 0; blue, 0 }  ,draw opacity=0.11 ]  ;
	\draw  [color={rgb, 255:red, 254; green, 0; blue, 0 }  ,draw opacity=1 ][fill={rgb, 255:red, 255; green, 0; blue, 0 }  ,fill opacity=1 ] (228.44,241.56) .. controls (228.44,240.7) and (229.14,240) .. (230,240) .. controls (230.86,240) and (231.56,240.7) .. (231.56,241.56) .. controls (231.56,242.42) and (230.86,243.11) .. (230,243.11) .. controls (229.14,243.11) and (228.44,242.42) .. (228.44,241.56) -- cycle ;
	\draw [color={rgb, 255:red, 65; green, 117; blue, 5 }  ,draw opacity=1 ][line width=1.5]    (95.4,151.89) -- (194,151) ;

	\draw [color={rgb, 255:red, 65; green, 117; blue, 5 }  ,draw opacity=1 ][line width=1.5]    (194,151) -- (229,190) ;

	\draw [color={rgb, 255:red, 245; green, 166; blue, 35 }  ,draw opacity=1 ][line width=1.5]    (95.4,181.89) -- (230,182) ;

	\draw [color={rgb, 255:red, 245; green, 166; blue, 35 }  ,draw opacity=0.71 ][line width=1.5]    (230,182) -- (230.6,240.6) ;

	\draw [color={rgb, 255:red, 65; green, 117; blue, 5 }  ,draw opacity=0.54 ][line width=1.5]    (229,190) -- (230,240) ;

	\draw [color={rgb, 255:red, 65; green, 117; blue, 5 }  ,draw opacity=0.77 ][line width=1.5]    (231.56,241.56) -- (317,241) ;

	\draw [color={rgb, 255:red, 245; green, 166; blue, 35 }  ,draw opacity=0.44 ][line width=1.5]    (230,241.56) -- (317,241) ;

	\draw [color={rgb, 255:red, 74; green, 144; blue, 226 }  ,draw opacity=1 ][line width=1.5]    (95.4,91.89) -- (230,92) ;

	\draw [color={rgb, 255:red, 74; green, 144; blue, 226 }  ,draw opacity=1 ][line width=1.5]    (230,92) -- (230.6,109.6) ;

	\draw [color={rgb, 255:red, 74; green, 144; blue, 226 }  ,draw opacity=1 ][line width=1.5]    (230.6,109.6) -- (306,129) ;

	\draw  [color={rgb, 255:red, 255; green, 0; blue, 0 }  ,draw opacity=0.46 ] (188.2,151) .. controls (188.2,147.96) and (190.8,145.5) .. (194,145.5) .. controls (197.2,145.5) and (199.8,147.96) .. (199.8,151) .. controls (199.8,154.04) and (197.2,156.5) .. (194,156.5) .. controls (190.8,156.5) and (188.2,154.04) .. (188.2,151) -- cycle ; \draw  [color={rgb, 255:red, 255; green, 0; blue, 0 }  ,draw opacity=0.46 ] (189.9,147.11) -- (198.1,154.89) ; \draw  [color={rgb, 255:red, 255; green, 0; blue, 0 }  ,draw opacity=0.46 ] (198.1,147.11) -- (189.9,154.89) ;
	\draw  [color={rgb, 255:red, 255; green, 0; blue, 0 }  ,draw opacity=0.46 ] (224.2,182) .. controls (224.2,178.96) and (226.8,176.5) .. (230,176.5) .. controls (233.2,176.5) and (235.8,178.96) .. (235.8,182) .. controls (235.8,185.04) and (233.2,187.5) .. (230,187.5) .. controls (226.8,187.5) and (224.2,185.04) .. (224.2,182) -- cycle ; \draw  [color={rgb, 255:red, 255; green, 0; blue, 0 }  ,draw opacity=0.46 ] (225.9,178.11) -- (234.1,185.89) ; \draw  [color={rgb, 255:red, 255; green, 0; blue, 0 }  ,draw opacity=0.46 ] (234.1,178.11) -- (225.9,185.89) ;
	\draw  [color={rgb, 255:red, 255; green, 0; blue, 0 }  ,draw opacity=0.46 ] (224.2,92) .. controls (224.2,88.96) and (226.8,86.5) .. (230,86.5) .. controls (233.2,86.5) and (235.8,88.96) .. (235.8,92) .. controls (235.8,95.04) and (233.2,97.5) .. (230,97.5) .. controls (226.8,97.5) and (224.2,95.04) .. (224.2,92) -- cycle ; \draw  [color={rgb, 255:red, 255; green, 0; blue, 0 }  ,draw opacity=0.46 ] (225.9,88.11) -- (234.1,95.89) ; \draw  [color={rgb, 255:red, 255; green, 0; blue, 0 }  ,draw opacity=0.46 ] (234.1,88.11) -- (225.9,95.89) ;
	\draw  [color={rgb, 255:red, 255; green, 0; blue, 0 }  ,draw opacity=0.46 ] (243.2,75.4) .. controls (243.2,72.36) and (245.8,69.9) .. (249,69.9) .. controls (252.2,69.9) and (254.8,72.36) .. (254.8,75.4) .. controls (254.8,78.44) and (252.2,80.9) .. (249,80.9) .. controls (245.8,80.9) and (243.2,78.44) .. (243.2,75.4) -- cycle ; \draw  [color={rgb, 255:red, 255; green, 0; blue, 0 }  ,draw opacity=0.46 ] (244.9,71.51) -- (253.1,79.29) ; \draw  [color={rgb, 255:red, 255; green, 0; blue, 0 }  ,draw opacity=0.46 ] (253.1,71.51) -- (244.9,79.29) ;
	
	\draw (337,238.89) node   {$\epsilon $};
	\draw (66,77.89) node   {$x( \epsilon )$};
	\draw (91,248.49) node [scale=0.7]  {$0$};
	\draw (230,254) node [scale=0.7,color={rgb, 255:red, 255; green, 0; blue, 0 }  ,opacity=1 ]  {$9$};
	\draw (293,74.9) node [scale=0.8,color={rgb, 255:red, 255; green, 0; blue, 0 }  ,opacity=1 ] [align=left] {Node failure};
	\draw (143,102.93) node [scale=0.8] [align=left] {\textcolor[rgb]{0.29,0.56,0.89}{Node 1}};
	\draw (141,141.93) node [scale=0.8] [align=left] {\textcolor[rgb]{0.25,0.46,0.02}{Node 2}};
	\draw (143,192.93) node [scale=0.8] [align=left] {\textcolor[rgb]{0.96,0.65,0.14}{Node 3}};
	\end{tikzpicture}}
	\caption{ }
	\label{lossdisc}
\end{figure}
In particular, Figure \ref{lossdisc} (a) shows how the loss function varies piece-wise linearly until $\epsilon=9$, where it undergoes the jump discontinuity of size $\Delta l(c^*)$. 
On the other hand, from Figure \ref{lossdisc} (b) we can notice that all nodes are solvent for $\epsilon < 6.5$ while for $\epsilon \approx 6.5$ node $2$ goes bankrupt as its outflow falls below $w_2=3$. As the shock magnitude increases, we reach the discontinuity point at $\epsilon=9$ where the network suffers a dramatic crisis as nodes $1$ and $3$ suddenly default. Notice in particular how node $3$ goes from fully solvent ($x_3=w_3$) to completely insolvent ($x_3=0$) as the shock crosses the critical threshold $\epsilon=9$.
 \end{example}
 \medskip
\subsection{Sensitivity of Nash equilibria in constrained quadratic network games}\label{sec:sensistivity-games}
In the literature, the constrained quadratic games introduced in Section \ref{sec:game} are often studied \cite{Bramoulle.Kranton:2016} with the matrix $P$  parameterized as $P^{(\delta)}=\delta G$ where $G$ is some fixed matrix encoding the network interconnections and $\delta >0$ is a parameter describing the strength of the network interaction among the agents. If we put $\delta^*=\rho(G)^{-1}$, we have that $\rho(\delta G)<1$ for $\delta <\delta^*$. While Proposition \ref{proposition:uniqueness-irreducible}  implies that, for every fixed $\delta <\delta^*$, the network equilibrium is unique and continuous in the exogenous flow $c$, its sensitivity to the variations of $c$ may grow unbounded when $\delta$ approaches $\delta^*$. As it turns out, this occurs when the limit network has multiple equilibria. Indeed, we have the following result showing that in this case, arbitrarily small variations in the exogenous flow $c$ will determine, for $\delta$ close to $\delta^*$, a variation in the equilibrium of the size of the set of equilibria for the limit case $\delta=\delta^*$.

\begin{corollary}
 For an irreducible matrix $G$ in $\R_+^{n\times n}$ and a vector $w$ in $\R_+^n$, and $\delta$ in $(0,\delta^*]$, where $\delta^*=1/\rho(G)^{-1}$, let $P^{(\delta)}=\delta G$ and let $\ov x^{(\delta)}(c)$ and $\ul x^{(\delta)}(c)$ to be the minimal and maximal network equilibrium of the network $(P^{(\delta)}, w)$ with exogenous flow $c$ in $\R^n$. Also, write $x^{(\delta)}$ for the network equilibrium when it is unique. Let $c^*$ be an exogenous flow  such that the $(P^{\delta^*},w)$ has multiple network equilibria. Then, 
\be\label{eq:uniformbound}\sup\limits_{\delta<\delta^*}\sup\limits_{c\,:\, \|c-c^*\|\leq\epsilon}\|x^{(\delta)}(c)-x^{(\delta)}(c^*)\|\geq \|\ov x^{(\delta^*)}(c^*)-\ul x^{(\delta^*)}(c^*)\|>0\,,\ee
for every monotone norm $\|\,\cdot\,\|$  and every $\eps>0$. 
\end{corollary}
\begin{proof}
It follows from the comparative statics in Proposition \ref{proposition:c} (iv) that, for $\delta<\delta^*$, 
\be\label{boundary1}x^{(\delta)}(c^*)\leq \ul x^{(\delta^*)}(c^*)\lneq \ov x^{(\delta^*)}(c^*)\,.\ee
Let $p$ be any left dominant eigenvector of $G$ and thus of all $P^{(\delta)}$. It then follows from Proposition \ref{proposition:uniqueness-irreducible} that, 
$$\ul x^{(\delta^*)}(c^*+\epsilon p)= \ov x^{(\delta^*)}(c^*+\epsilon p)\,,\qquad \forall\epsilon >0\,,$$ and thus, by Theorem \ref{theo:continuity} and Proposition \ref{proposition:c} (v) again,
\be\label{boundary2}\lim\limits_{\delta\downarrow\delta^*}x^{(\delta)}(c^*+\epsilon p)=\ov x^{(\delta^*)}(c^*+\epsilon p)\geq \ov x^{(\delta^*)}(c^*)\,.\ee
For every monotone norm  $\|\,\cdot\,\|$, \eqref{boundary1} and \eqref{boundary2} imply that
$$
\lim\limits_{\delta\downarrow\delta^*}\| \ul x^{(\delta)}(c^*+\epsilon p)-x^{(\delta)}(c^*)\|\geq \|\ov x^{(\delta^*)}(c^*)-\ul x^{(\delta^*)}(c^*)\|>0\,,
$$
so that \eqref{eq:uniformbound} holds true for every $\epsilon >0$. 
\qed\end{proof}

 \section{Conclusion}\label{conclusions}
This paper has analyzed network equilibria modeled as the solutions of a linear fixed point equation with saturation non-linearities. Necessary and sufficient conditions for uniqueness and a general expression describing all such equilibria for a general network with spectral radius not larger than $1$ have been proved. Finally, the dependence of the network equilibria on the exogenous flows in the network has been studied highlighting the existence of jump discontinuities. This model was first considered to determine clearing payments in the context of networked financial institutions interconnected by obligations and it is one of the simplest continuous model where shock propagation phenomena and cascading failure effects may occur. It also describes the Nash equilibria of constrained quadratic network games with strategic complementarities. Our results contribute to an in-depth analysis of such applications.

The understanding of the extent to which the network topology determines the structure of the solutions as well the possibility of these cascading effects to occur is still not sufficiently understood. As a future project, we aim at studying this for random networks with prescribed degree distributions.

\bibliographystyle{informs2014}
\bibliography{lit} 

\begin{APPENDIX}{Proof of Proposition \ref{prop:contraction}}\label{proof:contraction}
We start with the following result. 

\begin{lemma}\label{lemma:stable} Let $P$ in $\R_+^{n\times n}$ be a non-negative square matrix such that 
\begin{itemize}
\item there exists a non-negative vector $v\ne 0$ such that $Pv\leq v$;
\item for every $i=1,\ldots,n$, there exists a path in $\mc G_P$ connecting $i$ to some $j$ such that $(Pv)_j<v_j$.
\end{itemize}
Then, $\rho(P)<1$.
\end{lemma}
\begin{proof}
Notice that, for every $h\ge0$, $P^{h}v\le v$, so that, for $t\ge h$, non-negativity of $P^{t-h}$ implies that $(P^tv)=P^{t-h}P^{h}v\le P^{h}v$. 
On the other hand, existence of a length-$l_i$ path from $i$ to $j$  in $\mc G_P$ is equivalent to that $(P^{l_i})_{ij}>0$. 
Therefore, if there exists a length-$l_i$ path in $\mc G_P$ from $i$ to some $j$ such that $(Pv)_j<v_j$, then, for every $t> l_i$,
$$(P^tv)_i\le (P^{l_i+1}v)_i=\sum_{k=1}^n(P^{l_i})_{ik}(Pv)_k=\sum_{k=1}^n(P^{l_i})_{ik}v_k<(P^{l_i}v)_i\le v_i\,.$$

  Therefore, with $t=1+\max_il_i$, we have $(P^{t}v)_i<x_i$ for every $i$. Since $x_i>0$ for every $i$, we can find $\epsilon >0$ such that $P^tv\leq (1-\epsilon)v$. This implies that $\lim P^{ t m}= 0$ as $m$ grows large and thus $\rho(P^{t})<1$. This yields $\rho(P)<1$.\qed
\end{proof}\medskip

We can now proceed to the proof of Proposition  \ref{prop:contraction}. 

First, we prove existence of a positive vector $v$ satisfying \eqref{Pv<=v} for every non-expansive network. We proceed by induction on the number $s$ of classes of $P$. If $s=1$, i.e., $P$ if is irreducible, the result follows from  Proposition \ref{proposition:PF} (iii). Now, assume that the result holds true for $s-1$ and let us prove it for $s$. Consider the block structure \eqref{canonical} and notice that by the inductive  hypothesis we can find vectors $x^{(l)}$ of dimension $|\mc V_l|$ for $l=2,\dots ,s$  with all positive entries such that 
$$\sum_{h=l}^s P^{(lh)}v^{(h)}\leq v^{(l)}\,.$$ 
We now show that we can find $\alpha >0$ and  $x^{(1)}$ of dimension $|\mc V_1|$ with all positive entries, such that 
\be\label{eq:inequality}P^{(11)}v^{(1)}+\alpha\sum_{j= 2}^s P^{(1j)}v^{(j)}\leq v^{(1)}\,.\ee
Indeed, if $\rho(P^{(11)})<1$ this simply follows from a continuity argument. Instead, if $\rho(P^{(11)})=1$, 
then since $P^{(11)}$ is irreducible, it admits a positive right dominant eigenvalue $v^{(1)}=P^{(11)}v^{(1)}$ by Proposition \ref{proposition:PF} (iii). 
On the other hand, since $\mc V_1$ is final, we have that $P^{(1h)}=0$ for every $h=1,\ldots,s$, so that \eqref{eq:inequality} is satisfied as an equality for all possible values of $\alpha >0$. This implies that the vector $v=(v^{(1)},\dots ,v^{(s-1)}, v^{(s)})$ has all positive entries and satisfies $Pv\leq v$. 

Finally, we prove that, existence of a positive vector $v$ satisfying \eqref{Pv<=v}  implies that the network is non-expansive. From \eqref{Pv<=v}, using the fact that all entries of $v$ are strictly positive, we deduce that $P^t$ is a bounded sequence, so that $\rho(P)\leq 1$.
Now, assume that $\mc V_l$ is a non final class such that $\rho(P^{(ll)})=1$. Indicating as usual with $v^{(l)}$ the restriction of $v$ to $\mc V_l$, we obtain the relation
$$P^{(ll)}v^{(l)}+\sum_{h=l+1}^s P^{(lh)}v^{(h)}\leq v^{(l)}$$
from which we deduce that 
$P^{(ll)}v^{(l)}\lneq v^{(l)}$. 
Since $P^{(ll)}$ is irreducible, we can apply Lemma \ref{lemma:stable} and conclude that $\rho(P^{(ll)})<1$. 
\qed

\end{APPENDIX}

\section*{Acknowledgments.}
This work was partially supported by a MIUR grant ``Dipartimenti di Eccellenza 2018--2022'' [CUP: E11G18000350001], a MIUR  Research Project PRIN 2017 ``Advanced Network Control of Future Smart Grids'' (http://vectors.dieti.unina.it), the Swedish Research Council [2015-04066], and by the Compagnia di San Paolo.

\end{document}